\documentclass[11pt]{article}

\usepackage[margin=1in]{geometry}
\usepackage{amsmath,amssymb,amsthm,mathtools}
\usepackage{enumitem}
\usepackage{microtype}
\usepackage{graphicx}
\usepackage{hyperref}
\usepackage[T1]{fontenc}
\usepackage{lmodern}
\usepackage{xcolor}

\hypersetup{colorlinks=true,linkcolor=blue,citecolor=blue,urlcolor=blue}
\mathtoolsset{showonlyrefs}

\newtheorem{theorem}{Theorem}
\newtheorem{proposition}[theorem]{Proposition}
\newtheorem{lemma}[theorem]{Lemma}
\newtheorem{appendixlemma}{Lemma}[section]
\newtheorem{definition}[theorem]{Definition}

\theoremstyle{remark}

\newcommand{\E}{\mathbb E}
\newcommand{\Pp}{\mathbb P}
\newcommand{\R}{\mathbb R}
\newcommand{\1}{\mathbf 1}
\newcommand{\Tr}{\operatorname{Tr}}
\newcommand{\Var}{\operatorname{Var}}
\newcommand{\Comp}{\operatorname{Comp}}
\newcommand{\Cat}{\operatorname{Cat}}
\newcommand{\ip}[2]{\langle #1,#2\rangle}
\newcommand{\scm}{\mu_{\rm sc}}
\newcommand{\eps}{\varepsilon}
\newcommand{\proofstep}[1]{\par\medskip\noindent\textbf{#1}\par\smallskip}

\title{Spectra of sparse high-dimensional random geometric graphs}
\author{Yifan Cao\footnote{Department of Operations Research and Financial Engineering, Princeton University, yc0570@princeton.edu} \qquad  Yizhe Zhu\footnote{Department of Mathematics, University of Southern California, yizhezhu@usc.edu}}
\date{}

\begin{document}
\maketitle

\begin{abstract}
We determine the limiting empirical spectral distribution of sparse
high-dimensional random geometric graphs.  The vertices are independent
uniform points on the unit sphere \(S^{d-1}\), and two
vertices are joined when their inner product exceeds a threshold chosen to
give edge density \(p\).  The edges
therefore have the same marginal probabilities as in an Erd\H{o}s--R\'enyi
graph, but the latent geometry introduces dependence among them.  We show that
these correlations are asymptotically invisible to the global spectrum in two
sparse regimes.  If \(p\to0\), \(np\to\infty\), and
\(d=\Omega(np\log(1/p))\), then the empirical spectral distribution of
\(A/\sqrt{np}\) converges in probability to the semicircle law.  If
\(p=\alpha/n\) for a fixed \(\alpha>0\) and \(d=\omega(\log n)\), then the
empirical spectral distribution of \(A/\sqrt{\alpha}\) converges in probability
to the limiting spectral distribution of \(\mathcal G(n,\alpha/n)\).  The proof
combines the moment method with a cluster expansion that decomposes geometric
dependence into weak local interactions, allowing us to control every fixed
walk pattern in the moment calculation.
\end{abstract}

\section{Introduction}

Random geometric graphs arise naturally in models of networks with latent spatial
or geometric structure \cite{penrose2003random,liu2023probabilistic}, including
wireless communication networks \cite{mao2012connectivity}, biological networks
\cite{higham2008fitting}, and social graphs \cite{hoff2002latent}.  In these
models, each vertex is assigned a random position in a geometric space, and edges
are formed according to pairwise proximity.  This construction induces
dependencies between edges and distinguishes random geometric graphs from
classical Erd\H{o}s--R\'enyi graphs, where edges are independent.

Let \(X_1,\dots,X_n\in \R^d\) be independent random vectors uniformly distributed
on the unit sphere \(S^{d-1}\).  The threshold spherical random geometric graph
\(\mathcal G(n,d,p)\) is obtained by joining \(i\) and \(j\) when
\(\ip{X_i}{X_j}\ge \tau\).  The threshold \(\tau=\tau_{d,p}\) is chosen so that
\(\E \1_{\{\ip{X_i}{X_j}\ge \tau\}}=p\) for \(i\ne j\).
Thus \(\mathcal G(n,d,p)\) has the same edge density as the Erd\H{o}s--R\'enyi
graph \(\mathcal G(n,p)\), but its edges are generally dependent.

A central feature of the model is that the graph is generated from latent vertex
features.  This makes it useful in data-science settings where vertices are
represented by high-dimensional feature vectors and edges encode similarity
\cite{racz2017basic}.  The same geometry also creates the triangle-closing
property \cite{dall2002random}: two vertices that are both close to a common
neighbor are more likely to be close to each other.  Such correlations are absent
from Erd\H{o}s--R\'enyi graphs but are important for modeling clustering in
real-world networks.

Our goal is to determine the limiting spectral distribution of the
adjacency matrix of \(\mathcal G(n,d,p)\), and in particular to identify when it
agrees with that of \(\mathcal G(n,p)\).  Two features make the sparse
high-dimensional regime delicate: the edge indicators are nonlinear threshold
functions, and edges sharing latent points are dependent.

\paragraph{Testing high-dimensional geometry}
Testing whether an observed graph contains latent high-dimensional geometry
has motivated a substantial line of work on this model.  For fixed \(n\),
Devroye et al.~\cite{devroye2011high} showed that \(\mathcal G(n,d,p)\) and
\(\mathcal G(n,p)\) become indistinguishable in total variation as
\(d\to\infty\).  When \(n,d\to\infty\) and \(p\in(0,1)\) is fixed, Bubeck et
al.~\cite{bubeck2016testing} identified the sharp transition \(d\asymp n^3\).
They also initiated the sparse regime \(p=\alpha/n\), proving that the total
variation distance tends to \(1\) when \(d=o(\log^3 n)\) and conjecturing a
transition at \(d\asymp\log^3 n\).  Liu et al.~\cite{liu2022testing} proved the
complementary indistinguishability result under
\(d=\Omega(\log^{36}n)\).  More recently, Du et
al.~\cite{du2026resolution} established the conjectured detection threshold in
a broad regime with \(d>n\).

The testing problem has since been developed in several directions.  Baguley
et al.~\cite{baguley2025testing} analyzed random toroidal graphs using
Edgeworth-type expansions, and Brennan et al.~\cite{brennan2024threshold}
determined detection thresholds for anisotropic latent geometries.  Bangachev
and Bresler~\cite{bangachev2024fourier,bangachev2025sandwiching} studied
low-degree polynomial tests and sandwiching couplings.  Mao, Wu, and
Xu~\cite{mao2026smooth} obtained a sharp detection threshold for smooth
spherical kernels and formulated a related spectral conjecture.  On the
recovery side, Li and Schramm~\cite{li2023spectral} and Mao and
Zhang~\cite{mao2024impossibility} located the threshold for estimating a latent
inner product at \(d\asymp np\log(1/p)\).

Total variation concerns the entire graph distribution.  Here we ask the
coarser, but still structurally informative, question of whether the latent
geometry changes the global spectrum.
For the Erd\H{o}s--R\'enyi graph \(\mathcal G(n,p)\), the empirical spectral
distribution of the centered and normalized adjacency matrix converges to the
semicircle law in the sparse regime
\cite{furedi1981eigenvalues,bai2010spectral,tran2013sparse}.  If \(p\to0\) and
\(np\to\infty\), the uncentered matrix \(A/\sqrt{np}\) has the same limit:
centering changes it only by a rank-one matrix and a scalar shift that vanishes
as \(n\to\infty\).  This leads to the following question:
\begin{center}
    \textit{Can one distinguish \(\mathcal G(n,d,p)\) from \(\mathcal G(n,p)\)
    by their limiting spectral distributions?}
\end{center}

Our two main theorems answer this question in the negative, under their
respective dimension assumptions.  When \(p\to0\) and \(np\to\infty\), the
empirical spectral distribution of \(A/\sqrt{np}\) converges to the semicircle
law provided \(d=\Omega(np\log(1/p))\).  When \(p=\alpha/n\) for fixed
\(\alpha>0\), the empirical spectral distribution of \(A/\sqrt\alpha\)
converges to \(\nu_\alpha\), the limiting spectral distribution of
\(\mathcal G(n,\alpha/n)\), provided \(d=\omega(\log n)\).
The remaining case \(np\to0\) is degenerate: with high probability,
\((1-o(1))n\) vertices are isolated, so a \(1-o(1)\) fraction of the eigenvalues
are zero.  In the bounded expected degree regime, the Erd\H{o}s--R\'enyi limit
is not semicircular and has been studied from both physics and probability
perspectives
\cite{rodgers1988density,khorunzhy2004eigenvalue,bordenave2010resolvent}.

Thus Theorems~\ref{thm:main} and~\ref{thm:constantp} establish
universality at the level of the global empirical spectral distribution.
Motivated by the entropic and spectral thresholds
\(d\asymp np\log(1/p)\) in
\cite{bangachev2025sandwiching,li2023spectral,mao2024impossibility}, we
conjecture that the dimension scale in Theorem~\ref{thm:main} is sharp.

\paragraph{Nonlinear random matrices}
A nonlinear random matrix is obtained by applying a nonlinear function entrywise
to a classical random matrix ensemble.  In recent years, a rich theory has
emerged around such models, including kernel random matrices
\cite{bordenave2008eigenvalues,cheng2013spectrum,do2013spectrum,el2010spectrum,mei2022generalization,pandit2024universality},
as well as conjugate kernels and neural tangent kernels arising from neural
networks
\cite{louart2018random,pennington2017nonlinear,fan2020spectra,guionnet2025global}.
Semicircle laws have been obtained for broad classes of nonlinear random matrices
\cite{lu2022equivalence,dubova2023universality,wang2024deformed}.

Closest to our setting are the random inner-product kernel matrices studied in
\cite{lu2022equivalence,dubova2023universality}.  Let \(K\) be an \(n\times n\)
random matrix with
\(K_{ij}=n^{-1/2}f(\sqrt d\,\ip{X_i}{X_j})\1_{\{i\ne j\}}\), where
\(X_1,\dots,X_n\) are independent uniform points on \(S^{d-1}\), and \(f\) is
independent of \(d\).  In the regime \(d=\omega(n)\), Dubova et
al.~\cite{dubova2023universality} showed that the limiting spectral
distribution is a possibly rescaled semicircle law.  The normalized adjacency
matrix \(A/\sqrt n\) corresponds to the step kernel
\(f_{d,p}(x)=\1_{\{x\ge a_{d,p}\}}\), where
\(a_{d,p}=\sqrt d\,\tau_{d,p}\).  For fixed \(p\), this threshold converges to
a fixed Gaussian quantile.  In the sparse regime \(p\to0\), however,
\(a_{d,p}\to\infty\), so the kernel varies with \(n\) and the existing proofs
do not apply.

To the best of our knowledge, Theorems~\ref{thm:main}
and~\ref{thm:constantp} are the first results describing the limiting spectral
distributions of \(\mathcal G(n,d,p)\) in simultaneous high-dimensional
\((d\to\infty)\) and sparse \((p\to0)\) regimes.  The model combines two
difficult features of modern random matrix theory, nonlinearity and sparsity,
and requires correlation estimates that are uniform in the sparse regime.
Beyond the growing-dimension setting, the spectrum of fixed-dimensional random
geometric graphs has also been studied in
\cite{adhikari2022spectrum,hamidouche2023normalized}.

\section{Main results}\label{sec:main}

Let \(X_1,\ldots,X_n\) be independent uniform points on the unit sphere
\(S^{d-1}\subset \R^d\).  Let \(p=p_n\in(0,1/2)\), and let
\(\tau=\tau_{d,p}\) be the threshold satisfying
\[
        \Pp\{\ip{X_1}{X_2}\ge \tau\}=p .
\]
The spherical threshold random geometric graph has adjacency matrix
\[
        A_{ij}=\1_{\{\ip{X_i}{X_j}\ge \tau\}}
        \quad (i\ne j),\qquad
        A_{ii}=0 .
\]
Put \(q=1-p\), \(L=\log(1/p)\), and, for the centered part of the proof,
\[
        M_n=\frac{A-p(J-I)}{\sqrt{npq}} .
\]
Equivalently, with
\[
        \xi(x,y)=\1_{\{\ip{x}{y}\ge \tau\}}-p,
        \qquad
        \xi_{ij}=\xi(X_i,X_j),
\]
we have \((M_n)_{ij}=\xi_{ij}/\sqrt{npq}\) for \(i\ne j\), and
\((M_n)_{ii}=0\).

For an \(n\times n\) real symmetric matrix \(B\), write
\(\mu_B=n^{-1}\sum_{i=1}^n\delta_{\lambda_i(B)}\) for its empirical spectral
distribution.  The standard semicircle law \(\scm\) has density
\[
        \frac{1}{2\pi}\sqrt{4-x^2}\,\1_{\{|x|\le2\}}.
\]
All convergence statements below refer to weak convergence in probability of
these random probability measures.

\begin{theorem}[Sparse semicircle law]\label{thm:main}
Fix any constant \(c_0>0\).  Assume \(p\to0\),
\(np\to\infty\), and that, for all sufficiently large \(n\),
\[
        d\ge c_0 np\log(1/p).
\]
Then \(\mu_{A/\sqrt{np}}\) converges in probability to the semicircle law
\(\scm\).
\end{theorem}

We conjecture that the dimension condition in Theorem~\ref{thm:main} is optimal
for convergence to the semicircle law.  Determining the limiting spectral
distribution of \(A/\sqrt{np}\) when \(d=o(np\log(1/p))\) remains an
interesting open problem.

We next state the analogue in the bounded expected degree regime.  Here the
right comparison object is the sparse Erd\H{o}s--R\'enyi limiting law rather
than the semicircle law.

Let \(\nu_\alpha\) denote the limiting empirical spectral distribution of
\(A_{\mathrm{ER}}/\sqrt{\alpha}\), where \(A_{\mathrm{ER}}\) is the adjacency
matrix of the Erd\H{o}s--R\'enyi graph \(\mathcal G(n,\alpha/n)\)
\cite{bordenave2010resolvent,jung2018delocalization}.

\begin{theorem}[Bounded expected degree]\label{thm:constantp}
Let \(A\) be the adjacency matrix of \(\mathcal G(n,d,p)\).  Assume
\[
        p=\frac{\alpha}{n}
        \quad\text{for a fixed \(\alpha>0\)},\qquad
        d=\omega(\log n).
\]
Then \(\mu_{A/\sqrt{np}}=\mu_{A/\sqrt{\alpha}}\) converges in probability to
\(\nu_{\alpha}\).
\end{theorem}

\begin{figure}[!ht]
        \centering
        \includegraphics[width=\textwidth]{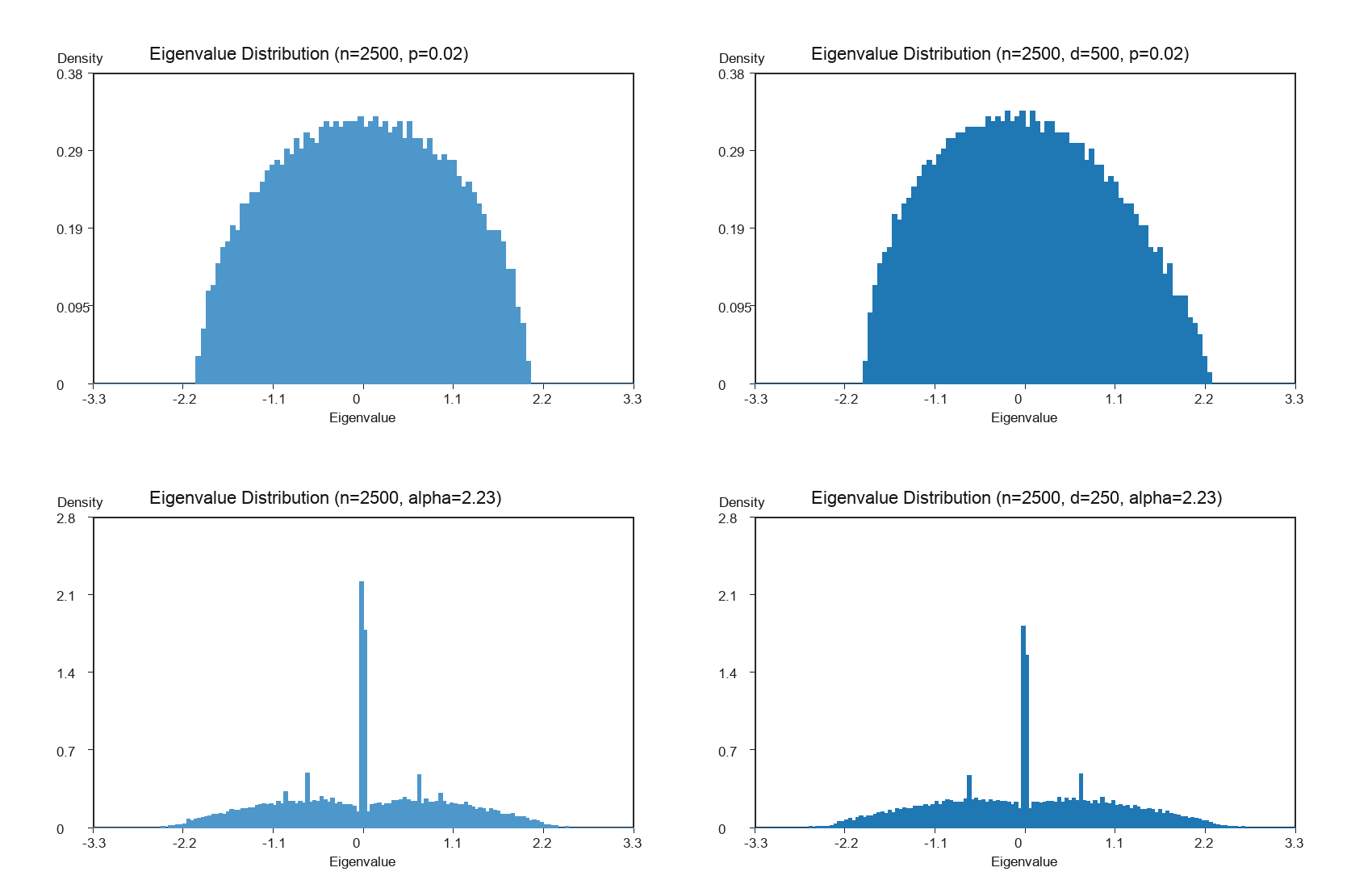}
        \caption{Simulations illustrating the two main regimes.  The
        histograms use \(n=2500\), and the thresholds for the geometric graphs
        are estimated from the one-dimensional spherical marginal.  Top:
        empirical spectral distributions of
        \((A-p(J-I))/\sqrt{np(1-p)}\) in the diverging expected degree regime
        for \(\mathcal G(n,p)\) and \(\mathcal G(n,d,p)\), with \(p=0.02\) and
        \(d=500\).  Bottom: empirical spectral distributions of
        \(A/\sqrt{\alpha}\) in the bounded expected degree regime for
        \(\mathcal G(n,\alpha/n)\) and \(\mathcal G(n,d,\alpha/n)\), with
        \(\alpha=2.23\) and \(d=250\).}
        \label{fig:main-simulations}
\end{figure}

Figure~\ref{fig:main-simulations} gives a finite-size comparison
of the geometric and Erd\H{o}s--R\'enyi spectra in the two regimes.  The
simulations are consistent with the agreement predicted by the theorems.

\paragraph{Technical overview}
The main obstacle is that the edge indicators exposed by a fixed walk may
be correlated through their shared latent points.  Proposition~\ref{prop:corr}
quantifies this dependence.  The centered correlation is zero whenever the
support graph has a leaf; when every component has minimum degree at least two,
it decays by a power of \(d^{-1}\) determined by the support graph.

The proof separates the geometric and combinatorial parts of this
estimate.  Lemma~\ref{lem:gram-density} compares the joint law of finitely many
spherical inner products with the product of their one-dimensional marginals.
A cluster expansion of the resulting density ratio
\cite{brydges1986short} expresses geometric dependence as products of local
generator monomials.  Informally, it decomposes the geometric dependence into
finite collections of weak local interactions.  Each interaction has a
dimension cost, and these costs add when interactions are combined.  Centering
and symmetry determine which collections may survive, and
Lemma~\ref{claim:cost} converts those constraints into the required power of
\(d^{-1}\).  Every surviving product contains a bounded admissible core, and
Lemma~\ref{lem:integrated-core} sums over all such cores.

Related high-trace arguments appear in
\cite{li2023spectral,liu2023local}.  Li and Schramm
\cite{li2023spectral} use a Gegenbauer expansion and orthogonality to contract
degree-two vertices.  Liu et al.~\cite{liu2023local} peel trees from the walk
support, contract degree-two paths in its \(2\)-core, and use the mixing rate of
the spherical-cap walk.  Their graph-correlation estimates suffice for the
operator-norm problems considered there, but do not provide bounds strong
enough to show that all non-Catalan contributions are negligible in our
global-law argument.

For a global law, one must isolate the exact Catalan contribution and
show that every other fixed walk support remains negligible after summing over
its vertex labels.  Proposition~\ref{prop:corr} provides the
support-sensitive factor \(d^{-\rho(H)}\) needed for this summation.  The rest
of the proof is a moment convergence argument: Wigner tree words give the semicircle
moments in the diverging expected degree regime, while tree-supported moments
agree exactly with the sparse Erd\H{o}s--R\'enyi moments in the bounded expected
degree regime.
Variance estimates then yield convergence in probability.

\paragraph{Organization of the paper}
Section~\ref{sec:correlation-bound} proves the finite-graph correlation
bound, which is the main technical input.  Sections~\ref{sec:proof_semicircle}
and~\ref{sec:proof_constant} then treat the diverging and bounded expected degree
regimes, respectively.  Appendix~\ref{sec:spherical-estimates} collects the
one-dimensional spherical estimates used in the correlation argument.

\paragraph{Notation}
Throughout the proofs, \(k\) denotes a fixed moment order.  Subscripts record
the dependence of constants on fixed parameters; for example, \(C_m\) depends
only on \(m\), \(C_k\) only on \(k\), and \(C_K,A_K,B_K\) only on \(K\).  When
a bound also depends on the fixed lower-bound constant \(c_0\) or the fixed
parameter \(\alpha\), we write \(C_{k,c_0}\) or \(C_{k,\alpha}\), respectively.
Constants with the same subscript may change from line to line, but never
depend on \(n,p,d\).  Absolute constants are written as \(c_*,C_*\) when their
role is reused.  The \(r\)-th Catalan number is denoted by
\(\Cat_r=\frac{1}{r+1}\binom{2r}{r}\), \(r\ge0\).  Since \(p\to0\), we use
\(q=1-p\ge1/2\) for all sufficiently large \(n\).

\section{Finite-graph correlation bound}\label{sec:correlation-bound}

Throughout this section, \(H\) is a fixed simple graph.  We define
\(V(H)\) to be the set of endpoints of its edges, thereby discarding isolated
vertices.  When \(H\) is nonempty, relabel these vertices so that
\(V(H)=\{1,\ldots,m\}\); for the empty graph, set \(V(H)=\varnothing\) and
\(m=0\).  Set
\(a=\sqrt d\,\tau\) and \(Z_{ij}=\sqrt d\,\ip{X_i}{X_j}\).  For
\(e=\{i,j\}\), write \(Z_e=Z_{ij}\) and
\(\xi_e(z)=\1_{\{z\ge a\}}-p\).  Thus
\(\xi(X_i,X_j)=\xi_{\{i,j\}}(Z_{ij})\).

The following proposition is the central estimate of this section.  A
leaf produces an exact cancellation; when every component has minimum degree at
least two, the remaining dependence decays as an explicit power of \(d^{-1}\).

For such a graph \(H\), define its centered edge correlation by
\(\Phi(H):=\E\prod_{\{u,v\}\in E(H)}\xi(X_u,X_v)\).
\begin{proposition}[Finite centered correlation estimate]\label{prop:corr}
Fix an integer \(K\ge1\), and let \(L=\log(1/p)\) and
\(\eps_n=\sqrt{L/d}\).  Assume that \(p\to0\) and \(L/d\to0\).  There are
constants \(A_K,B_K<\infty\), depending only on \(K\), such that every simple
graph \(H\) with at most \(K\) edges satisfies the following.
\begin{enumerate}
    \item If \(H\) has a leaf, then \(\Phi(H)=0\).
    \item If every connected component of \(H\) has minimum degree at least
    \(2\), then
\begin{align}
        |\Phi(H)|
        \le
        A_K p^{|E(H)|-B_K\eps_n}L^{B_K}d^{-\rho(H)}.
        \label{eq:correlation-bound}
\end{align}
Here \(\rho(H)\) is specified by
\eqref{eq:rho-connected}--\eqref{eq:rho-general}: for a connected graph \(C\),
set
\begin{align}
        \rho(C)=\max\left\{\frac12,\ |V(C)|-1-\frac{|E(C)|}{2}\right\}.
       \label{eq:rho-connected}
\end{align}
For an arbitrary graph \(H\), set
\begin{align}
        \rho(H)=\sum_{C\in\Comp(H)}\rho(C),
        \label{eq:rho-general}
\end{align}
where the sum is over the connected components of \(H\).
\end{enumerate}
\end{proposition}

We prove the two assertions of Proposition~\ref{prop:corr} separately.
The leaf cancellation follows directly from conditioning.  For the second
assertion, we first reduce to connected \(H\) and then analyze the Gram-density
ratio.

\begin{proof}[Proof of the first assertion in Proposition~\ref{prop:corr}]
Suppose that \(v\) is a leaf of \(H\), with unique neighbor \(u\).
Condition on \((X_w)_{w\ne v}\).  All factors in \(\Phi(H)\) are then fixed
except \(\xi(X_v,X_u)\).  For every fixed \(x\in S^{d-1}\), rotational
invariance of the uniform probability measure \(\sigma\) on \(S^{d-1}\) gives
\[
        \int_{S^{d-1}}\xi(y,x)\,d\sigma(y)
        =0.
\]
Therefore the conditional expectation of the full product is zero, and
hence \(\Phi(H)=0\).
\end{proof}

It remains to prove \eqref{eq:correlation-bound}.  If \(H\) is empty,
then \(\Phi(H)=1\) and \(\rho(H)=0\), so the estimate follows after taking
\(A_K\ge1\).  Now let \(H\) be nonempty, with connected components
\(C_1,\ldots,C_s\).  Their vertex sets are disjoint, and hence independence of
the underlying points and \eqref{eq:rho-general} give
\(\Phi(H)=\prod_{j=1}^s\Phi(C_j)\) and
\(\rho(H)=\sum_{j=1}^s\rho(C_j)\).

It therefore suffices to prove \eqref{eq:correlation-bound} for connected
\(H\).  Indeed, multiplying the component bounds replaces the prefactor
\(A_K\) by \(A_K^s\) and the loss exponent \(B_K\) by \(sB_K\); since
\(s\le K\), both changes are absorbed by enlarging \(A_K,B_K\).  Henceforth
assume that \(H\) is connected and has minimum degree at least two.  With
\(h=|E(H)|\), the degree-sum identity gives
\(2h=\sum_{v\in V(H)}\deg_H(v)\ge2m\), so \(m\le h\le K\).

Let \(\nu_d\) be the one-dimensional law of \(Z_{12}\), and let
\(R_{m,d}\) be the density ratio of the joint law of
\((Z_{ij})_{1\le i<j\le m}\) with respect to
\(\nu_d^{\otimes \binom m2}\).  By a change of measure,
\begin{align}
        \Phi(H)
        =\int \prod_{e\in E(H)}\xi_e(Z_e)R_{m,d}(Z)
        \,d\nu_d^{\otimes \binom m2}(Z).
        \label{eq:Phi-integral}
\end{align}

The proof uses three main ingredients.  Lemma~\ref{lem:gram-density} identifies the
density ratio that rewrites the correlation as a product-measure integral.
Lemmas~\ref{lem:log-density-ratio} and~\ref{lem:cluster} expand that ratio into
generator products, after which centering and parity determine which products
can survive.  Finally, Lemma~\ref{claim:cost} supplies the required dimension
gain, while the absolute-sum bound in Lemma~\ref{lem:cluster} and the integrated
summation in Lemma~\ref{lem:integrated-core} control all products containing
the resulting bounded cores.

\subsection{The Gram determinant and the log-determinant expansion}

We begin with an explicit density formula for finitely many spherical inner
products.  This formula is standard in multivariate analysis
\cite{muirhead1982aspects}, but the proof is included for
completeness.
\begin{lemma}[Joint density of finitely many spherical Gram variables]
\label{lem:gram-density}
Let \(X_1,\ldots,X_m\) be independent uniform random points on
\(S^{d-1}\), and assume \(d\ge m\).  Let \(B=B(Z)\) be the \(m\times m\)
symmetric matrix with zero diagonal and off-diagonal entries \(B_{ij}=Z_{ij}\).
Then the joint density \(f_{m,d}\) of the off-diagonal variables
        \((Z_{ij})_{1\le i<j\le m}\)
satisfies
\begin{align}
        f_{m,d}(Z)
        =
        C_{m,d}\,
        \det(I_m+d^{-1/2}B(Z))^{(d-m-1)/2}
        \1_{\{I_m+d^{-1/2}B(Z)\succ0\}},
        \label{eq:joint_density_Z_correct}
\end{align}
where, with \(N=\binom m2\),
\begin{equation}
        C_{m,d}
        =
        d^{-N/2}
        \prod_{j=2}^m
        \frac{\Gamma(d/2)}
        {\pi^{(j-1)/2}\Gamma((d-j+1)/2)}.
        \label{eq:joint-density-constant}
\end{equation}
\end{lemma}

\begin{proof}
Define the ordinary Gram matrix
\begin{align}
        G=(G_{ij})_{1\le i,j\le m},
        \qquad
        G_{ij}=\ip{X_i}{X_j}.
\end{align}
Then
\begin{align}
        G=I_m+d^{-1/2}B(Z).
        \label{eq:G-from-Z}
\end{align}
It is enough to prove that the off-diagonal entries of \(G\) have density
proportional to
       \(\det(G)^{(d-m-1)/2}\)
on the region \(G\succ0\), \(G_{ii}=1\).  Indeed, the change of variables
\(G_{ij}=d^{-1/2}Z_{ij}\), \(i<j\), has constant Jacobian
\(d^{-\binom m2/2}\).  Multiplying the conditional
normalizing constants in \eqref{eq:projection-density} and this Jacobian gives
\eqref{eq:joint-density-constant}.  Substituting \eqref{eq:G-from-Z} then gives
\eqref{eq:joint_density_Z_correct}.

We now derive the density of the off-diagonal entries of \(G\).  Expose
the vectors sequentially.  For \(1\le j\le m\), let
\(G_j=(\ip{X_a}{X_b})_{1\le a,b\le j}\) be the \(j\times j\) Gram
matrix of \(X_1,\ldots,X_j\).  Since \(d\ge m\),
the vectors \(X_1,\ldots,X_j\) are linearly independent almost surely for
every \(j\le m\), so \(G_j\) is positive definite almost surely.

Fix \(j\in\{2,\ldots,m\}\), and condition on
        \(X_1,\ldots,X_{j-1}\).
Let \(S_{j-1}=\operatorname{span}(X_1,\ldots,X_{j-1})\).  Choose an
orthonormal basis \(e_1,\ldots,e_{j-1}\) of \(S_{j-1}\).  Define the projection coordinate vector
\[
        Y_j
        =
        \begin{pmatrix}
        \ip{e_1}{X_j}\\
        \vdots\\
        \ip{e_{j-1}}{X_j}
        \end{pmatrix}
        \in\mathbb R^{j-1}.
\]
Because \(X_j\) is uniform on \(S^{d-1}\) and independent of the previously
exposed vectors, the projection \(Y_j\) has density on the unit ball
\(\{y\in\mathbb R^{j-1}:\|y\|<1\}\) given by
\begin{equation}
        C_{d,j}\,(1-\|y\|^2)^{(d-j-1)/2},
        \qquad
        C_{d,j}=
        \frac{\Gamma(d/2)}
        {\pi^{(j-1)/2}\Gamma((d-j+1)/2)}.
        \label{eq:projection-density}
\end{equation}
Formula~\eqref{eq:projection-density} follows from the coarea formula applied
to the coordinate projection \(S^{d-1}\to\mathbb R^{j-1}\); see
\cite[Chapter~3]{evans2015measure}.
Now define the vector of new inner products
\[
        s_j
        =
        \begin{pmatrix}
        \ip{X_1}{X_j}\\
        \vdots\\
        \ip{X_{j-1}}{X_j}
        \end{pmatrix}
        \in\mathbb R^{j-1}.
\]
Write the previously exposed vectors \(X_1,\ldots,X_{j-1}\) in the basis
\(e_1,\ldots,e_{j-1}\).  Let \(A_{j-1}\) be the \((j-1)\times(j-1)\)
matrix whose \(i\)-th column is the coordinate vector of \(X_i\).  Then
\begin{equation}
        A_{j-1}^{\top}A_{j-1}=G_{j-1}.
        \label{eq:A-gram}
\end{equation}
Also,
\begin{equation}
        s_j=A_{j-1}^{\top}Y_j.
        \label{eq:s-from-Y}
\end{equation}
Therefore \(Y_j=(A_{j-1}^{\top})^{-1}s_j\).
The Jacobian of the map \(Y_j\mapsto s_j=A_{j-1}^{\top}Y_j\) is
\[
        |\det A_{j-1}^{\top}|
        =
        |\det A_{j-1}|
        =
        \det(G_{j-1})^{1/2},
\]
where the last equality follows from \eqref{eq:A-gram}.  Hence
\begin{equation}
        dY_j
        =
        \det(G_{j-1})^{-1/2}\,ds_j.
        \label{eq:jacobian}
\end{equation}
Furthermore,
\begin{equation}
        \|Y_j\|^2
        =
        s_j^\top G_{j-1}^{-1}s_j.
        \label{eq:norm-Y}
\end{equation}
Combining \eqref{eq:projection-density}, \eqref{eq:jacobian}, and
\eqref{eq:norm-Y}, the conditional density of \(s_j\) given
\(G_{j-1}\) is
\begin{equation}
        C_{d,j}\,
        \det(G_{j-1})^{-1/2}
        \left(1-s_j^\top G_{j-1}^{-1}s_j\right)^{(d-j-1)/2},
        \label{eq:conditional-density-s}
\end{equation}
on the region \(s_j^\top G_{j-1}^{-1}s_j<1\).
The new Gram matrix is
\[
        G_j=
        \begin{pmatrix}
        G_{j-1} & s_j\\
        s_j^\top & 1
        \end{pmatrix}.
\]
By the Schur complement formula,
\begin{equation}
        \det(G_j)
        =
        \det(G_{j-1})
        \left(1-s_j^\top G_{j-1}^{-1}s_j\right).
        \label{eq:schur-complement}
\end{equation}
Thus \(1-s_j^\top G_{j-1}^{-1}s_j=\det(G_j)/\det(G_{j-1})\).
Substituting this into \eqref{eq:conditional-density-s}, we obtain
\begin{equation}
\begin{aligned}
        f(s_j\mid G_{j-1})
        &=
        C_{d,j}\,
        \det(G_{j-1})^{-1/2}
        \left(
        \frac{\det(G_j)}{\det(G_{j-1})}
        \right)^{(d-j-1)/2}  \\
        &=
        C_{d,j}\,
        \det(G_j)^{(d-j-1)/2}
        \det(G_{j-1})^{-(d-j)/2}.
\end{aligned}
        \label{eq:conditional-density-det}
\end{equation}

The concatenated vector \((s_2,\ldots,s_m)\) is exactly the vector of
off-diagonal entries of \(G_m\), up to a permutation of coordinates.  Hence
this change of ordering has Jacobian one.  Moreover, \(G_{j-1}\) is determined
by \((s_2,\ldots,s_{j-1})\), so the chain rule for conditional densities and
\eqref{eq:conditional-density-det} give
\begin{equation}
\begin{aligned}
 f_{(s_2,\ldots,s_m)}(s_2,\ldots,s_m)
 &=\prod_{j=2}^m f(s_j\mid s_2,\ldots,s_{j-1}) \\
 &=\prod_{j=2}^m f(s_j\mid G_{j-1}) \\
 &=\left(\prod_{j=2}^m C_{d,j}\right)
   \prod_{j=2}^m
   \det(G_j)^{(d-j-1)/2}
   \det(G_{j-1})^{-(d-j)/2}
   \1_{\{G_m\succ0\}}.
\end{aligned}
        \label{eq:telescope-product}
\end{equation}
The support conditions \(G_j\succ0\), \(2\le j\le m\), are equivalent to
\(G_m\succ0\), since the \(G_j\) are its leading principal submatrices.  The
determinant product in \eqref{eq:telescope-product} telescopes:
\[
\begin{aligned}
 &\prod_{j=2}^m
   \det(G_j)^{(d-j-1)/2}
   \det(G_{j-1})^{-(d-j)/2} \\
 &\quad=
   \det(G_m)^{(d-m-1)/2}\det(G_1)^{-(d-2)/2}
   \prod_{j=2}^{m-1}
   \det(G_j)^{(d-j-1)/2-(d-j-1)/2} \\
 &\quad=\det(G_m)^{(d-m-1)/2},
\end{aligned}
\]
where the last identity uses \(G_1=[1]\).  Thus the off-diagonal entries of
\(G=G_m\) have density proportional to
\(\det(G)^{(d-m-1)/2}\1_{\{G\succ0\}}\).
This proves Lemma~\ref{lem:gram-density}.
\end{proof}

The joint density in Lemma~\ref{lem:gram-density} is useful because it can be
compared directly with the product of the one-edge marginals.
Lemma~\ref{lem:log-density-ratio} records this comparison in logarithmic
form.
\begin{lemma}[Log-density ratio]
\label{lem:log-density-ratio}
On the positive definite region
        \(I_m+d^{-1/2}B\succ0\),
one has
\begin{equation}
\begin{aligned}
        \log R_{m,d}(Z)
        &=
        \frac{d-m-1}{2}
        \log\det(I_m+d^{-1/2}B)
        -
        \frac{d-3}{2}
        \sum_{i<j}
        \log\left(1-\frac{Z_{ij}^2}{d}\right)
        +
        \kappa_{m,d},
\end{aligned}
        \label{eq:log-density-ratio}
\end{equation}
where \(\kappa_{m,d}\) is a constant depending only on \(m\) and \(d\).
\end{lemma}

\begin{proof}
Let \(N=\binom m2\).  The product density is
\[
\begin{aligned}
        f_{\mathrm{prod}}(Z)
        &=
        \prod_{i<j}
        c_d
        \left(1-\frac{Z_{ij}^2}{d}\right)^{(d-3)/2}
        =
        c_d^N
        \prod_{i<j}
        \left(1-\frac{Z_{ij}^2}{d}\right)^{(d-3)/2}
        .
\end{aligned}
\]
On the region \(I_m+d^{-1/2}B\succ0\), all \(2\times2\) principal minors are
positive.  Since the diagonal of
\(I_m+d^{-1/2}B\) is \(1\), this implies
        \(1-\frac{Z_{ij}^2}{d}>0\)
for every \(i<j\).  Hence \(|Z_{ij}|<\sqrt d\), so the product marginal density
is strictly positive there.
By Lemma~\ref{lem:gram-density}, specifically
\eqref{eq:joint_density_Z_correct}, on the positive definite region,
\[
\begin{aligned}
        R_{m,d}(Z)
        &=
        \frac{
        C_{m,d}\,
        \det(I_m+d^{-1/2}B)^{(d-m-1)/2}
        }{
        c_d^N
        \prod_{i<j}
        \left(1-\frac{Z_{ij}^2}{d}\right)^{(d-3)/2}
        } .
\end{aligned}
\]
Taking logarithms gives
\[
\begin{aligned}
        \log R_{m,d}(Z)
        &=
        \log C_{m,d}
        -
        N\log c_d
        +
        \frac{d-m-1}{2}
        \log\det(I_m+d^{-1/2}B)
        -
        \frac{d-3}{2}
        \sum_{i<j}
        \log\left(1-\frac{Z_{ij}^2}{d}\right).
\end{aligned}
\]
Define
        \(\kappa_{m,d}
        =
        \log C_{m,d}
        -
        N\log c_d\).
This proves the claim.
\end{proof}

\subsection{Integrated cluster expansion}
\label{subsec:integrated-cluster}

We now expand the exact density ratio into monomials, track their powers of
\(d^{-1}\), and establish the bounds needed to integrate the resulting series.

Lemma~\ref{lem:log-density-ratio} separates \(\log R_{m,d}\) into a determinant
term and a sum of scalar, one-coordinate terms.  When
\(\|d^{-1/2}B\|<1\),
\[
 \log\det(I_m+d^{-1/2}B)
 =\Tr\log(I_m+d^{-1/2}B)
 =\sum_{\ell\ge1}\frac{(-1)^{\ell+1}}{\ell}
   d^{-\ell/2}\Tr(B^\ell).
\]
Since \(B_{ii}=0\), one has \(\Tr B=0\), and
\[
 \Tr(B^\ell)
 =\sum_{i_0,\ldots,i_{\ell-1}=1}^m
   B_{i_0i_1}B_{i_1i_2}\cdots B_{i_{\ell-1}i_0}.
\]
Every nonzero summand is therefore the monomial of a closed walk of length
\(\ell\).  After multiplication by the prefactor \((d-m-1)/2\) in
\eqref{eq:log-density-ratio}, the terms with \(\ell\ge3\) have bounded
coefficients times \(d^{1-\ell/2}\).  For example, a triangle contributes a
multiple of \(d^{-1/2}Z_{12}Z_{23}Z_{31}\).

For each pair \(i<j\), the scalar term has the expansion
\[
 -\frac{d-3}{2}\log\left(1-\frac{Z_{ij}^2}{d}\right)
 =\sum_{s\ge1}\frac{d-3}{2ds}\,d^{1-s}Z_{ij}^{2s}.
\]
The terms with \(s\ge2\) have bounded coefficients times
\(d^{1-s}Z_{ij}^{2s}\).  The \(s=1\) term combines with the \(\ell=2\)
determinant term according to
\[
 -\frac{d-m-1}{2d}\sum_{i<j} Z_{ij}^2
 +\frac{d-3}{2d}\sum_{i<j} Z_{ij}^2
 =\frac{m-2}{2d}\sum_{i<j} Z_{ij}^2.
\]
Consequently, apart from the constant \(\kappa_{m,d}\), the logarithmic
density ratio consists of closed-walk monomials with \(\ell\ge3\), quadratic
monomials, and higher even powers of one coordinate.  Up to bounded
coefficients (uniformly in \(d\) for fixed \(m\)), their respective powers of
\(d^{-1}\) are
\[
        d^{-(\ell/2-1)},\qquad d^{-1},\qquad d^{-(s-1)}.
\]
The next paragraph formalizes these three families and defines the cost of each
term to be the corresponding exponent of \(d^{-1}\).

\paragraph{Generators}
Fix \(m\ge 2\), and write \(\mathcal E_m=\{\{i,j\}:1\le i<j\le m\}\).  For
\(e=\{i,j\}\in\mathcal E_m\), write \(Z_e=Z_{ij}\).
A \emph{closed-walk generator} is a closed walk
\(w=(i_0,i_1,\ldots,i_{\ell-1},i_\ell=i_0)\) in the complete graph on
\(\{1,\ldots,m\}\), with length \(\ell\ge3\), such that
\(i_{r-1}\ne i_r\) for every \(r\).  Its monomial is
\(M_w(Z)=Z_{i_0i_1}Z_{i_1i_2}\cdots Z_{i_{\ell-1}i_0}\), and its cost is
\(\omega(w)=\ell/2-1\).

A \emph{one-edge correction generator} is a pair \((e,s)\), where
\(e\in\mathcal E_m\) and \(s\ge1\) is an integer.  Its monomial is
\(M_{e,s}(Z)=Z_e^{2s}\), and its cost is
\[
        \omega(e,s)=
        \begin{cases}
        1, & s=1,\\
        s-1, & s\ge2.
        \end{cases}
\]

A \emph{product of generators} is a finite multiset
\(\Gamma=\{g_1,\ldots,g_r\}\), where each \(g_a\) is either a closed-walk
generator or a one-edge correction generator.  We define
\(Z^\Gamma=\prod_{a=1}^r M_{g_a}(Z)\) and
\(\omega(\Gamma)=\sum_{a=1}^r \omega(g_a)\).
The empty product \(\Gamma=\varnothing\) is allowed; then
\(Z^\varnothing=1\) and \(\omega(\varnothing)=0\).

If \(\Gamma_0\) and \(\Gamma\) are products of generators, we say that
\(\Gamma\) \emph{contains} \(\Gamma_0\) as a subproduct if every generator
appearing in \(\Gamma_0\) appears in \(\Gamma\) with at least the same
multiplicity.

With these definitions, write the nonconstant part of \(\log R_{m,d}\) as
\[
        U(Z)=\sum_g\lambda_gd^{-\omega(g)}M_g(Z).
\]
On the region where this series is absolutely convergent, as proved in
Lemma~\ref{lem:cluster},
\[
\begin{aligned}
        e^{U(Z)}
        &=
        \sum_{r=0}^{\infty}\frac1{r!}
        \left(\sum_g\lambda_gd^{-\omega(g)}M_g(Z)\right)^r\\
        &=
        \sum_{r=0}^{\infty}\frac1{r!}
        \sum_{g_1,\ldots,g_r}
        \left(\prod_{j=1}^r\lambda_{g_j}\right)
        d^{-\sum_{j=1}^r\omega(g_j)}
        \prod_{j=1}^rM_{g_j}(Z).
\end{aligned}
\]
For \(\Gamma=\{g_1,\ldots,g_r\}\), the definitions give
\[
        \prod_{j=1}^rM_{g_j}(Z)=Z^\Gamma,
        \qquad
        \sum_{j=1}^r\omega(g_j)=\omega(\Gamma).
\]
Thus every such product carries the factor
\(d^{-\omega(\Gamma)}Z^\Gamma\): costs add when monomials are multiplied.
This bookkeeping is analogous to a cluster expansion
\cite{brydges1986short}, in which each generator represents a local interaction
among the inner-product coordinates.

After establishing the analytic bounds for this expansion, we will
combine it with centering and symmetry.  The exact product-law factorization
leading to the coverage and parity conditions is given immediately before
Definition~\ref{def:admissible}.
Lemma~\ref{lem:cluster} gives the precise expansion and integrability
estimates used in \eqref{eq:Phi-integral}.
\begin{lemma}[Integrated Gram expansion]\label{lem:cluster}
Fix integers \(m\ge2\) and \(N\ge0\).  There are constants \(C_m<\infty\),
\(C_{m,N}<\infty\), and \(\delta_m>0\) such that, for all sufficiently large
\(d\), the first two statements below hold for every fixed
\(\delta\in(0,\delta_m]\).
\begin{enumerate}
    \item On the coordinatewise event
\[
        \mathcal G_\delta
        =\left\{\max_{i<j}|Z_{ij}|\le\delta\sqrt d\right\},
\]
the density ratio has an absolutely convergent expansion
\begin{equation}
        R_{m,d}(Z)
        =
        \sum_{\Gamma}
        c_\Gamma d^{-\omega(\Gamma)}Z^\Gamma.
        \label{eq:R-expansion}
\end{equation}
The sum is over all finite products of the generators defined above, including
the empty product.  The coefficients \(c_\Gamma\) may depend on \(m\) and \(d\).
\item For \(Z\in\mathcal G_\delta\), if \(\Gamma_0\) contains at most \(N\) generators, counted with
multiplicity, then the absolute
sum of all terms in \eqref{eq:R-expansion} containing \(\Gamma_0\) as a
subproduct is bounded by
\begin{equation}
        C_{m,N} d^{-\omega(\Gamma_0)}|Z^{\Gamma_0}|
        \exp\left(
        C_m\sum_{i<j}\frac{|Z_{ij}|^3}{\sqrt d}
        \right).
        \label{eq:absolute-subproduct-bound}
\end{equation}
\item Fix \(\delta\in(0,\delta_m]\).  Let \(p=p_n\in(0,1)\), set
\(L=\log(1/p)\), and assume
\(p\to0\), \(d\to\infty\), and \(L/d\to0\).  Choose \(a\) so that
\(\Pp\{Z_{12}\ge a\}=p\).
Let \(\mu\) be either the true joint law of \((Z_{ij})_{i<j}\) or the product
law \(\nu_d^{\otimes\binom m2}\).  For every fixed \(A,\Lambda>0\) and
integers \(K,D\ge0\), there is a constant
\(C_{A,K,m,D,\Lambda,\delta}<\infty\) such that, for all sufficiently large
\(n\), every set \(F\subseteq\mathcal E_m\) with \(|F|\le K\) and every
polynomial \(W\) of degree at most \(D\) and coefficient \(\ell^1\)-norm at most
\(\Lambda\) satisfy
\[
        \left|
        \int_{\mathcal G_\delta^c}
        W(Z)\prod_{e\in F}\xi_e(Z_e)\,d\mu(Z)
        \right|
        \le C_{A,K,m,D,\Lambda,\delta}p^Ad^{-A},
\]
where \(\xi_e(z)=\1_{\{z\ge a\}}-p\).
\end{enumerate}
\end{lemma}

\begin{proof}
\leavevmode\par
\proofstep{Step 1: Expansion of the log-density ratio.}
By Lemma~\ref{lem:log-density-ratio}, specifically
\eqref{eq:log-density-ratio}, it is enough to expand the determinant term and
the one-coordinate logarithms on a region where both series converge.
Set \(\delta_m=(4m^2)^{-1}\) and
\(\delta=\delta_m\) throughout Steps~1--5.  Then
\((m-1)\delta<1/2\).  Since \(B\) is symmetric, its operator norm is at most
its maximum absolute row sum.  Hence, on \(\mathcal G_\delta\),
\(\|d^{-1/2}B\|\le(m-1)\delta<1/2\).
Every eigenvalue of \(I_m+d^{-1/2}B\) is therefore at least
\(1-\|d^{-1/2}B\|>1/2\), so this matrix is positive definite.  Moreover,
\[
\begin{aligned}
        \sum_{\ell\ge1}\frac{d^{-\ell/2}}{\ell}
        |\Tr(B^\ell)|
        &\le
        m\sum_{\ell\ge1}\frac{\|d^{-1/2}B\|^\ell}{\ell}
        <\infty.
\end{aligned}
\]
Thus the logarithmic series is absolutely convergent and gives
\[
        \log\det(I+d^{-1/2}B)
        =
        \sum_{\ell\ge1}
        \frac{(-1)^{\ell+1}}{\ell}d^{-\ell/2}\Tr(B^\ell),
\]
where the equality follows by applying the scalar logarithmic series to the
eigenvalues of \(d^{-1/2}B\).  Since \(B_{ii}=0\),
\[
        \Tr B=0,
        \qquad
        \Tr B^2
        =\sum_{i,j}B_{ij}B_{ji}
        =2\sum_{i<j}Z_{ij}^2.
\]
The determinant part gives
\begin{equation}
\begin{aligned}
        \frac{d-m-1}{2}\log\det(I_m+d^{-1/2}B)
        &=
        \frac{d-m-1}{2}
        \sum_{\ell\ge2}
        \frac{(-1)^{\ell+1}}{\ell}
        d^{-\ell/2}\Tr B^\ell  \\
        &=
        -\frac{d-m-1}{4d}\Tr B^2
        +
        \sum_{\ell\ge3}
        a_{\ell,d}\,d^{1-\ell/2}\Tr B^\ell,
\end{aligned}
        \label{eq:det-expansion}
\end{equation}
where, for \(\ell\ge3\),
\[
        a_{\ell,d}
        =\frac{d-m-1}{2d}\frac{(-1)^{\ell+1}}{\ell},
        \qquad
        |a_{\ell,d}|\le\frac1\ell
\]
for all sufficiently large \(d\), with \(m\) fixed.  For the scalar terms, use
\[
        -\log(1-u)=\sum_{s\ge1}\frac{u^s}{s},
        \qquad |u|<1.
\]
On \(\mathcal G_\delta\), every \(|Z_{ij}|/\sqrt d\le \delta\), so this series
is absolutely convergent because
\[
        \sum_{s\ge1}\frac1s
        \left|\frac{Z_{ij}^2}{d}\right|^s
        \le
        \sum_{s\ge1}\frac{\delta^{2s}}s
        <\infty.
\]
Consequently,
\begin{equation}
\begin{aligned}
        -\frac{d-3}{2}
        \sum_{i<j}
        \log\left(1-\frac{Z_{ij}^2}{d}\right)
        &=
        \frac{d-3}{2}
        \sum_{i<j}
        \sum_{s\ge1}
        \frac{1}{s}
        \left(\frac{Z_{ij}^2}{d}\right)^s  \\
        &=
        \frac{d-3}{2d}
        \sum_{i<j}Z_{ij}^2
        +
        \sum_{i<j}\sum_{s\ge2}
        b_{s,d}\,d^{1-s}Z_{ij}^{2s},
\end{aligned}
        \label{eq:marginal-expansion}
\end{equation}
where, for \(s\ge2\),
\[
        b_{s,d}=\frac{d-3}{2ds},
        \qquad
        |b_{s,d}|\le\frac1{2s}
\]
for all sufficiently large \(d\).  Using
\(\Tr B^2=2\sum_{i<j}Z_{ij}^2\), the quadratic contributions from
\eqref{eq:det-expansion} and \eqref{eq:marginal-expansion}, respectively, are
\[
        -\frac{d-m-1}{2d}\sum_{i<j}Z_{ij}^2,
        \qquad
        \frac{d-3}{2d}\sum_{i<j}Z_{ij}^2.
\]
Their sum is
\begin{equation}
        \frac{m-2}{2d}\sum_{i<j}Z_{ij}^2.
        \label{eq:quad-correction}
\end{equation}
\proofstep{Step 2: Identification of the generators and their costs.}
Recall the cost definition at the beginning of
Section~\ref{subsec:integrated-cluster}: it is the exponent of \(d^{-1}\) in
the coefficient of a generator.  The coefficient in
\eqref{eq:quad-correction} is \((m-2)d^{-1}/2\), so each summand is a
one-edge correction generator of cost \(1\).  For \(\ell\ge3\),
\[
        \Tr B^\ell
        =
        \sum_{i_0,\ldots,i_{\ell-1}=1}^m
        B_{i_0i_1}B_{i_1i_2}\cdots B_{i_{\ell-1}i_0}.
\]
Because \(B_{ii}=0\), a nonzero summand has
\(i_{r-1}\ne i_r\) for every \(r\), with indices read cyclically.  It is
therefore the monomial
\(Z_{i_0i_1}\cdots Z_{i_{\ell-1}i_0}\) of a closed walk of length \(\ell\),
and its prefactor in \eqref{eq:det-expansion} is exactly
\[
        a_{\ell,d}d^{1-\ell/2}
        =
        a_{\ell,d}d^{-(\ell/2-1)}.
\]
Thus it is a closed-walk generator of cost \(\ell/2-1\).

The remaining one-coordinate terms in \eqref{eq:marginal-expansion} have
the form \(d^{1-s}Z_{ij}^{2s}\), \(s\ge2\).
These are one-edge correction generators.  Their cost is
       \(s-1\ge1\).
Together with the quadratic correction \eqref{eq:quad-correction}, these are
exactly the one-edge correction generators from the definition.

Hence the nonconstant part of \(\log R_{m,d}\) may be written as
\begin{equation}
        U(Z)
        =
        \sum_{g}
        \lambda_g d^{-\omega(g)}M_g(Z),
        \label{eq:U-generator-expansion}
\end{equation}
where \(g\) ranges over all generators, \(M_g\) is its monomial,
and \(\omega(g)\) is its cost.  Thus
\[
        \log R_{m,d}(Z)=\kappa_{m,d}+U(Z),
        \qquad
        R_{m,d}(Z)=e^{\kappa_{m,d}}e^{U(Z)}.
\]
More explicitly, \(\lambda_g=a_{\ell,d}\) for a
closed walk of length \(\ell\), \(\lambda_g=(m-2)/2\) for a quadratic
one-edge correction, and \(\lambda_g=b_{s,d}\) for a one-edge correction with
\(s\ge2\).  Thus \(|\lambda_g|\le C_m\).  The constant term
\(\kappa_{m,d}\) is kept separately.  Every closed-walk generator
has length \(\ell\ge3\), and therefore
\[
 \omega(g)=\frac{\ell}{2}-1\ge\frac12.
\]
For a one-edge correction, \(\omega(e,1)=1\), while
\(\omega(e,s)=s-1\ge1\) for \(s\ge2\).  Hence every generator satisfies
\(\omega(g)\ge1/2\).

\proofstep{Step 3: Absolute convergence of the generator series.}
For the absolute convergence needed below, set
\[
        \mathcal A(Z)
        =
        \sum_g |\lambda_g|d^{-\omega(g)}|M_g(Z)|.
\]
We first bound the closed-walk contribution to \(\mathcal A\).  Let
\(N_m=\binom m2\), \(S=(\sum_{i<j}Z_{ij}^2)^{1/2}\), and let
\(\mathcal W_\ell\) be the set of closed-walk generators of length \(\ell\).
Choosing the \(\ell\) successive vertices gives
\[
        |\mathcal W_\ell|\le m^\ell,
        \qquad
        |M_w(Z)|\le S^\ell\quad (w\in\mathcal W_\ell).
\]
Using \(|\lambda_w|\le1/\ell\le1\), we work on \(\mathcal G_\delta\), where
\(S/\sqrt d\le N_m^{1/2}\delta\).  Since
\(N_m^{1/2}\le m\), our choice of \(\delta_m\) gives
\(mN_m^{1/2}\delta\le1/4\).  The geometric
series bound and the norm comparison
\(S^3\le N_m^{1/2}\sum_{i<j}|Z_{ij}|^3\) then give
\[
\begin{aligned}
        \sum_{\ell\ge3}\sum_{w\in\mathcal W_\ell}
        |\lambda_w|d^{1-\ell/2}|M_w(Z)|
        &\le d\sum_{\ell\ge3}
        \left(\frac{mS}{\sqrt d}\right)^\ell \\
        &\le \frac{2m^3S^3}{\sqrt d}
        \le C_m\sum_{i<j}\frac{|Z_{ij}|^3}{\sqrt d}.
\end{aligned}
\]

Next consider the one-edge correction terms.  The quadratic correction gives
\(C_m\sum_{i<j}Z_{ij}^2/d\).  The elementary bound
\(t^2/d\le d^{-1}+t^3/\sqrt d\), valid for \(t\ge0\) and \(d\ge1\), gives
\[
        C_m\sum_{i<j}\frac{Z_{ij}^2}{d}
        \le
        \frac{C_mN_m}{d}
        +
        C_m\sum_{i<j}\frac{|Z_{ij}|^3}{\sqrt d}.
\]
For the higher one-edge corrections, using \(|Z_{ij}|\le \delta\sqrt d\),
\[
        d^{1-s}|Z_{ij}|^{2s}
        =
        \frac{|Z_{ij}|^3}{\sqrt d}
        \left(\frac{|Z_{ij}|}{\sqrt d}\right)^{2s-3}
        \le
        \frac{|Z_{ij}|^3}{\sqrt d}\delta^{2s-3}.
\]
Therefore, for each edge \(e=\{i,j\}\),
\[
\begin{aligned}
        \sum_{s\ge2}|b_{s,d}|d^{1-s}|Z_{ij}|^{2s}
        &\le
        \frac{|Z_{ij}|^3}{\sqrt d}
        \sum_{s\ge2}\frac{\delta^{2s-3}}{2s} \\
        &\le
        C_\delta\frac{|Z_{ij}|^3}{\sqrt d}.
\end{aligned}
\]
Since \(N_m\), \(C_{\delta_m}\), and \(\delta_m\) depend only on \(m\), the
closed-walk and one-edge estimates, together with \(d^{-1}\le1\), give
\begin{equation}
        \mathcal A(Z)
        \le C_m\left(1+
        \sum_{i<j}\frac{|Z_{ij}|^3}{\sqrt d}\right).
        \label{eq:generator-sum-bound}
\end{equation}
For \(Z\in\mathcal G_\delta\),
\[
 \sum_{i<j}\frac{|Z_{ij}|^3}{\sqrt d}
 \le N_m\delta^3d<\infty.
\]
Thus \(\mathcal A(Z)<\infty\), which is exactly the absolute convergence of
\eqref{eq:U-generator-expansion}; moreover, \(|U(Z)|\le\mathcal A(Z)\).

\proofstep{Step 4: Normalization, exponentiation, and the first assertion.}
At \(Z=0\), the density-ratio formula gives
\[
        e^{\kappa_{m,d}}=\frac{C_{m,d}}{c_d^{N_m}}.
\]
Here
\(c_d=\Gamma(d/2)/(\sqrt{\pi d}\,\Gamma((d-1)/2))\).
Stirling's formula gives, uniformly for \(2\le j\le m\),
\[
        \frac{\Gamma(d/2)}{\Gamma((d-j+1)/2)}
        =
        \left(\frac d2\right)^{(j-1)/2}
        \left(1+O_m(d^{-1})\right).
\]
Substitution into the normalizing-constant formula
\eqref{eq:joint-density-constant} from Lemma~\ref{lem:gram-density}, together
with \(\sum_{j=2}^m(j-1)=N_m\), gives
\[
\begin{aligned}
        C_{m,d}
        &=(2\pi)^{-N_m/2}\left(1+O_m(d^{-1})\right),\\
        c_d^{N_m}
        &=(2\pi)^{-N_m/2}\left(1+O_m(d^{-1})\right),
\end{aligned}
\]
where the second line uses the case \(j=2\) of the same gamma-ratio estimate.
Consequently, \(e^{\kappa_{m,d}}=1+O_m(d^{-1})\), so
\(e^{\kappa_{m,d}}\le C_m\) for all sufficiently large \(d\).

Put \(u_g=\lambda_gd^{-\omega(g)}M_g(Z)\), and let \(\mathfrak N\) be the
set of nonnegative integer sequences \(\mathbf n=(n_g)_g\) with finite
support.  Since \(\sum_g|u_g|=\mathcal A(Z)<\infty\), Tonelli's theorem,
applied first to finite sets of generators and then by monotone convergence,
gives
\[
\begin{aligned}
        \sum_{\mathbf n\in\mathfrak N}
        \prod_g\frac{|u_g|^{n_g}}{n_g!}
        &=\prod_g\left(\sum_{n\ge0}\frac{|u_g|^n}{n!}\right)\\
        &=\exp\left(\sum_g|u_g|\right)
        =e^{\mathcal A(Z)}<\infty.
\end{aligned}
\]
Thus the multivariate exponential series is absolutely convergent.  Passing
to the limit in the corresponding finite product identity and identifying
\(\mathbf n\) with its generator multiset \(\Gamma\) gives
\[
        R_{m,d}(Z)
        =
        \sum_{\Gamma}
        \left(
        e^{\kappa_{m,d}}
        \prod_g\frac{\lambda_g^{\,n_g(\Gamma)}}{n_g(\Gamma)!}
        \right)
        d^{-\omega(\Gamma)}Z^\Gamma.
\]
Here \(n_g(\Gamma)\) is the multiplicity of \(g\) in \(\Gamma\), and only
finitely many factors differ from \(1\).  Therefore
\eqref{eq:R-expansion} holds with
\begin{equation}
        c_\Gamma
        =
        e^{\kappa_{m,d}}
        \prod_g\frac{\lambda_g^{\,n_g(\Gamma)}}{n_g(\Gamma)!}.
        \label{eq:cGamma-definition}
\end{equation}
This proves the first assertion.

\proofstep{Step 5: Terms containing a prescribed subproduct.}
We now prove \eqref{eq:absolute-subproduct-bound}.  Fix a product of generators
\(\Gamma_0\), and let \(n_g^0=n_g(\Gamma_0)\).  Set
\(x_g=|\lambda_g|d^{-\omega(g)}|M_g(Z)|\).  Since
\(\sum_gx_g=\mathcal A(Z)<\infty\), the factorial inequality
\((n_0+k)!\ge n_0!k!\) gives, for \(x\ge0\),
\[
        \sum_{n\ge n_0}\frac{x^n}{n!}
        \le \frac{x^{n_0}}{n_0!}e^x.
\]
The condition \(\Gamma\supseteq\Gamma_0\) is equivalent to
\(n_g(\Gamma)\ge n_g^0\) for every \(g\).  Hence
\eqref{eq:cGamma-definition}, Tonelli's theorem, and the preceding inequality
give
\[
\begin{aligned}
        \sum_{\Gamma\supseteq\Gamma_0}
        |c_\Gamma|d^{-\omega(\Gamma)}|Z^\Gamma|
        &=e^{\kappa_{m,d}}
        \prod_g\left(\sum_{n\ge n_g^0}\frac{x_g^n}{n!}\right) \\
        &\le e^{\kappa_{m,d}+\mathcal A(Z)}
        \left(\prod_g\frac{x_g^{n_g^0}}{n_g^0!}\right).
\end{aligned}
\]
If \(\Gamma_0\) contains at most \(N\) generators, then
\[
\begin{aligned}
        \prod_g\frac{x_g^{n_g^0}}{n_g^0!}
        &=
        d^{-\omega(\Gamma_0)}|Z^{\Gamma_0}|
        \prod_g\frac{|\lambda_g|^{n_g^0}}{n_g^0!}\\
        &\le
        C_{m,N}d^{-\omega(\Gamma_0)}|Z^{\Gamma_0}|.
\end{aligned}
\]
Here the identity uses the additive definitions of \(\omega(\Gamma_0)\) and
\(Z^{\Gamma_0}\), and the inequality uses \(|\lambda_g|\le C_m\) and
\(\sum_gn_g^0\le N\).  Finally, \eqref{eq:generator-sum-bound} and the boundedness of
\(e^{\kappa_{m,d}}\) yield
\[
        \sum_{\Gamma\supseteq\Gamma_0}
        |c_\Gamma|d^{-\omega(\Gamma)}|Z^\Gamma|
        \le
        C_{m,N}d^{-\omega(\Gamma_0)}|Z^{\Gamma_0}|
        \exp\left(
        C_m\sum_{i<j}\frac{|Z_{ij}|^3}{\sqrt d}
        \right).
\]
This is \eqref{eq:absolute-subproduct-bound}.

Steps~1--5 prove the first two assertions on
\(\mathcal G_{\delta_m}\).  If \(0<\delta\le\delta_m\), then
\(\mathcal G_\delta\subseteq\mathcal G_{\delta_m}\); restricting the same
pointwise expansion and bound to \(\mathcal G_\delta\) proves the first two
assertions for every such \(\delta\).

\proofstep{Step 6: The complement of the convergence region.}
It remains to prove the third assertion.  Fix
\(\delta\in(0,\delta_m]\).  Under either choice of \(\mu\), each
coordinate has marginal law \(\nu_d\).  Hence the union bound, symmetry, and
the Gaussian tail estimate \eqref{eq:basic-gaussian-tail} in
Lemma~\ref{lem:tail} give
\[
\begin{aligned}
        \mu(\mathcal G_\delta^c)
        &\le \binom m2\nu_d\{|Z_{12}|>\delta\sqrt d\}\\
        &\le C_m e^{-c_*\delta^2 d}.
\end{aligned}
\]
Every coordinate satisfies \(|Z_e|\le\sqrt d\) under both laws.  Thus every
monomial of degree at most \(D\) has absolute value at most \(d^{D/2}\), and
the coefficient assumption gives \(|W(Z)|\le\Lambda d^{D/2}\).  Since
\(|\xi_e|\le1\),
\begin{equation*}
\begin{aligned}
        \left|
        \int_{\mathcal G_\delta^c}
        W(Z)\prod_{e\in F}\xi_e(Z_e)\,d\mu(Z)
        \right|
        &\le \Lambda d^{D/2}\mu(\mathcal G_\delta^c)\\
        &\le C_m\Lambda d^{D/2}e^{-c_*\delta^2d}.
\end{aligned}
\end{equation*}
Finally, since \(L=o(d)\) and \(\log d=o(d)\),
\[
        \frac{C_m\Lambda d^{D/2}e^{-c_*\delta^2d}}{p^Ad^{-A}}
        =C_m\Lambda\exp\left(-c_*\delta^2d+AL+\left(A+\frac D2\right)\log d\right)
        \longrightarrow0.
\]
The preceding integral is therefore at most
\(C_{A,K,m,D,\Lambda,\delta}p^Ad^{-A}\) for all sufficiently large \(n\).
This proves the third assertion and completes the proof.
\end{proof}

\subsection{Admissibility and the covering cost}

We now isolate the combinatorial mechanism behind the dimension gain.  Let
\(V\) be a finite vertex set containing \(V(H)\), and write the monomial of a
generator product as
\(Z^\Gamma=\prod_{e\in E(K_V)}Z_e^{r_e(\Gamma)}\), where
\(r_e(\Gamma)\in\mathbb Z_{\ge0}\).  Under the full product law,
independence gives the exact factorization
\[
\begin{aligned}
        &\int Z^\Gamma\prod_{e\in E(H)}\xi_e(Z_e)
        \prod_{f\in E(K_V)}d\nu_d(Z_f)\\
        &\qquad=
        \prod_{e\in E(H)}
        \int \xi_e(z)z^{r_e(\Gamma)}\,d\nu_d(z)
        \prod_{e\in E(K_V)\setminus E(H)}
        \int z^{r_e(\Gamma)}\,d\nu_d(z).
\end{aligned}
\]
If \(e\in E(H)\) and \(r_e(\Gamma)=0\), then the corresponding factor is
\(\int\xi_e\,d\nu_d=p-p=0\).  Thus a nonzero integral requires
\(r_e(\Gamma)\ge1\) for every \(e\in E(H)\).  If
\(e\in E(K_V)\setminus E(H)\) and \(r_e(\Gamma)\) is odd, then symmetry of
\(\nu_d\) gives \(\int z^{r_e(\Gamma)}\,d\nu_d(z)=0\).  Hence a nonzero term
must cover \(H\) and have even exponent on every auxiliary edge.  We record
these two necessary conditions in the following definition.
\begin{definition}[Admissibility]\label{def:admissible}
With the notation above, \(\Gamma\) \emph{covers} \(H\) if the first condition
below holds, and it is \emph{\(H\)-admissible} if both conditions hold:
\[
        r_e(\Gamma)\ge1\quad(e\in E(H)),
        \qquad
        r_e(\Gamma)\equiv0\pmod 2
        \quad(e\in E(K_V)\setminus E(H)).
\]
\end{definition}

The absolute-sum bound in Lemma~\ref{lem:cluster} will be used after a
bounded subproduct has been singled out.  Lemma~\ref{lem:integrated-core}
supplies the uniform summation over all possible choices of that subproduct,
including generators of arbitrarily large degree.

Fix integers \(m\ge2\) and \(K,N\ge0\), a constant \(B>0\), and a simple graph \(H\) on a
subset of \(\{1,\ldots,m\}\), with \(h=|E(H)|\le K\).  For \(s\in[0,K]\), let
\(\mathfrak C_{N,s}(H)\) be the set of \(H\)-admissible generator
products \(\Gamma_0\) that contain at most \(N\) generators and satisfy
\(\omega(\Gamma_0)\ge s\).

For the one-coordinate estimates, let \(X,Y\) be independent uniform points on
\(S^{d-1}\), set \(Z=\sqrt d\,\ip{X}{Y}\), and denote its law by \(\nu_d\).
Choose \(a\) so that \(\Pp\{Z\ge a\}=p\), and set
\(L=\log(1/p)\), \(\xi(z)=\1_{\{z\ge a\}}-p\), and
\(\eps_n=\sqrt{L/d}\).

\begin{lemma}[Integrated core summation]\label{lem:integrated-core}
Assume that \(p\to0\) and \(L/d\to0\).  There exist
\(\delta=\delta_{m,B}\in(0,\delta_m]\), \(C=C_{m,B}<\infty\), and
\(C_{m,K,N,B}<\infty\) such that, for all sufficiently large \(n\) and every
integer \(r\ge0\),
\begin{align}
        &\E_{Z\sim\nu_d}\left[|\xi(Z)|\,|Z|^r
        e^{B|Z|^3/\sqrt d}\1_{\{|Z|\le\delta\sqrt d\}}\right]
        \le C^{r+1}p^{1-C\eps_n}(\sqrt L+\sqrt r)^r,
        \label{eq:weighted-uniform-centered}\\
        &\E_{Z\sim\nu_d}\left[|Z|^r
        e^{B|Z|^3/\sqrt d}\1_{\{|Z|\le\delta\sqrt d\}}\right]
        \le C^{r+1}(1+\sqrt r)^r.
        \label{eq:weighted-uniform-full}
\end{align}
Moreover,
\begin{equation}
\begin{aligned}
        &\sum_{\Gamma_0\in\mathfrak C_{N,s}(H)}
        d^{-\omega(\Gamma_0)}
        \int_{\mathcal G_\delta}
        \prod_{e\in E(H)}|\xi_e(Z_e)|\,
        |Z^{\Gamma_0}|
        \exp\left(B\sum_{e\in\mathcal E_m}
        \frac{|Z_e|^3}{\sqrt d}\right)
        \,d\nu_d^{\otimes\binom m2}(Z)  \\
        &\hspace{35mm}\le
        C_{m,K,N,B}\,
        p^{h-C_{m,K,N,B}\eps_n}L^{C_{m,K,N,B}}d^{-s}.
\end{aligned}
        \label{eq:integrated-core-sum}
\end{equation}
\end{lemma}

\begin{proof}
All constants below may depend on the fixed parameters indicated in the
statement.

\proofstep{Part 1: The uncentered one-coordinate bound.}
By the Gaussian density estimate \eqref{eq:basic-gaussian-density} in
Lemma~\ref{lem:tail}, there are absolute constants \(c_*,C_* >0\) such that
\(f_d(z)\le C_*e^{-c_*z^2}\) for \(|z|\le\sqrt d/2\).  Choose
\(\delta\in(0,\min\{\delta_m,1/2\}]\) so that
\(B\delta\le c_*/3\); we may decrease \(\delta\) further in Part~3.  On
\(|z|\le\delta\sqrt d\),
\[
        \frac{B|z|^3}{\sqrt d}\le B\delta z^2\le\frac{c_*}{3}z^2.
\]
It follows that
\[
\begin{aligned}
 &\E\left[|Z|^r e^{B|Z|^3/\sqrt d}
        \1_{\{|Z|\le\delta\sqrt d\}}\right] \\
 &\qquad\le C_*
        \int_{\mathbb R}|z|^r e^{-2c_*z^2/3}\,dz
 \le C^{r+1}(1+\sqrt r)^r,
\end{aligned}
\]
where the last step is the gamma-integral bound.  This proves
\eqref{eq:weighted-uniform-full}.

\proofstep{Part 2: The centered one-coordinate bound.}
Since \(|\xi(z)|\le\1_{\{z\ge a\}}+p\), the term carrying the factor \(p\) is
bounded by Part~1.  For the cap term, Lemma~\ref{lem:tail}, specifically
\eqref{eq:density-decay-from-a}, gives
\[
        f_d(z)\le C(1+a)p e^{-c_*(z^2-a^2)},
        \qquad z\ge a.
\]
Moreover, for \(a\le z\le\delta\sqrt d\),
\[
 \frac{B(z^3-a^3)}{\sqrt d}
 \le \frac{3Bz}{2\sqrt d}(z^2-a^2)
 \le \frac{c_*}{2}(z^2-a^2).
\]
Thus the shifted Gaussian integral estimate
\eqref{eq:shifted-gaussian-integral} in Lemma~\ref{lem:tail}, with the decay
constant decreased if necessary, yields
\[
\begin{aligned}
 &\int_a^{\delta\sqrt d}z^r e^{Bz^3/\sqrt d}f_d(z)\,dz \\
 &\qquad\le C(1+a)p e^{Ba^3/\sqrt d}
        \int_a^\infty z^r e^{-c_*(z^2-a^2)/2}\,dz \\
 &\qquad\le C^{r+1}p e^{Ba^3/\sqrt d}(a+\sqrt r)^r.
\end{aligned}
\]
By \eqref{eq:a-bound} in Lemma~\ref{lem:tail}, \(a\le C\sqrt L\), and hence
\[
        e^{Ba^3/\sqrt d}
        \le e^{CL\sqrt{L/d}}=p^{-C\eps_n}.
\]
Since \(L\to\infty\), the cap term and the term carrying the factor \(p\) are
together at most
\(C^{r+1}p^{1-C\eps_n}(\sqrt L+\sqrt r)^r\).  This proves
\eqref{eq:weighted-uniform-centered}.

\proofstep{Part 3: Summation over the cores.}
For a generator \(g\), write \(j(g)=2\omega(g)\), and let \(D(g)\) be the
degree of its monomial.  Every \(j(g)\) is a positive integer.  The definitions
of the two types of generators give
\[
        j(g)\le D(g)\le j(g)+2.
\]
Consequently, if \(\Gamma_0\) has \(q\le N\) generators, total degree \(D\),
and \(j=2\omega(\Gamma_0)\), then
\begin{equation}
        j\le D\le j+2q\le j+2N.
        \label{eq:core-degree-cost}
\end{equation}

We also need a uniform count of products of a given cost.  For each integer
\(t\ge1\), there are at most \(A_m^{t+1}\) generators \(g\) with
\(j(g)=t\): a closed-walk generator has length \(t+2\), and there are only
\(O_m(1)\) one-edge correction generators of a given cost.  An ordered list
with total cost \(j/2\) determines a composition
\(j=t_1+\cdots+t_q\) into positive integers.  Summing over all compositions,
and thereby overcounting the products with at most \(N\) generators, gives
\[
 \sum_{q=1}^j\ \sum_{\substack{t_1+\cdots+t_q=j\\t_i\ge1}}
        \prod_{i=1}^q A_m^{t_i+1}
 \le 2^{j-1}A_m^{2j}.
\]
Thus the number of generator products with \(2\omega(\Gamma_0)=j\) is at most
\(Q_m^{j+1}\) for a constant \(Q_m\) depending only on \(m\); for \(j=0\),
the empty product is the only possibility.

Let \(r_e\) denote the exponent of \(Z_e\) in \(Z^{\Gamma_0}\).  Since both
\(\mathcal G_\delta\) and the exponential weight factor over the coordinates,
the summand associated with \(\Gamma_0\) equals
\[
\begin{aligned}
 I(\Gamma_0)
 &:={d^{-j/2}}
 \int_{\mathcal G_\delta}
        \prod_{e\in E(H)}|\xi_e(Z_e)|\,|Z^{\Gamma_0}|
        e^{B\sum_e|Z_e|^3/\sqrt d}\,d\nu_d^{\otimes\binom m2}(Z) \\
 &=d^{-j/2}
 \prod_{e\in E(H)}
 \E\left[|\xi(Z)|\,|Z|^{r_e}e^{B|Z|^3/\sqrt d}
        \1_{\{|Z|\le\delta\sqrt d\}}\right] \\
 &\qquad\times
 \prod_{e\in\mathcal E_m\setminus E(H)}
 \E\left[|Z|^{r_e}e^{B|Z|^3/\sqrt d}
        \1_{\{|Z|\le\delta\sqrt d\}}\right].
\end{aligned}
\]
Applying Parts~1 and~2 and using \(r_e\le D\le j+2N\), we obtain
\[
\begin{aligned}
 I(\Gamma_0)
 &\le C_{m,N,B}Q_{m,B}^j p^{h-C_{m,K,B}\eps_n}d^{-j/2}
        (L+j)^{(j+2N)/2}.
\end{aligned}
\]
Indeed, the product of the moment factors is bounded by
\[
 \prod_{e\in E(H)}(\sqrt L+\sqrt{r_e})^{r_e}
 \prod_{e\notin E(H)}(1+\sqrt{r_e})^{r_e}
 \le C^D(L+j)^{D/2},
\]
and the \(h\le K\) centered coordinates contribute
\(p^{h-C_{m,K,B}\eps_n}\).  After including the number of products of cost
\(j/2\) and enlarging \(Q_{m,B}\), the total contribution at that cost is at
most
\[
 C_{m,K,N,B}Q_{m,B}^j p^{h-C_{m,K,N,B}\eps_n}
 d^{-j/2}(L+j)^{(j+2N)/2}.
\]

We first sum over \(2s\le j\le cd\), where \(c>0\) will be chosen below.  Since
\(s\le K\),
\[
\begin{aligned}
 Q_{m,B}^j d^{-j/2}(L+j)^{(j+2N)/2}
 &\le C_{m,K,B}d^{-s}(L+j)^{N+s}
 \left(\frac{Q_{m,B}^2(L+j)}d\right)^{j/2-s}.
\end{aligned}
\]
For \(j\le2L\), the expression in parentheses is at most
\(3Q_{m,B}^2L/d=o(1)\), so summing the resulting geometric series gives
\[
 \sum_{2s\le j\le2L}(L+j)^{N+s}
 \left(\frac{Q_{m,B}^2(L+j)}d\right)^{j/2-s}
 \le C_{m,K,N,B}L^{N+K}.
\]
If \(2L<j\le cd\), then \(L+j\le3j/2\).  Choose \(c=c_{m,B}>0\) so small
that \(3Q_{m,B}^2c/2\le\theta<1\).  Then
\[
 \sum_{2L<j\le cd}(L+j)^{N+s}
 \left(\frac{Q_{m,B}^2(L+j)}d\right)^{j/2-s}
 \le C\sum_{j\ge1}j^{N+K}\theta^{j/2-K}<\infty.
\]
The contribution of all \(j\le cd\) is therefore at most
\[
        C_{m,K,N,B}p^{h-C_{m,K,N,B}\eps_n}
        L^{C_{m,K,N,B}}d^{-s}.
\]

It remains to consider \(j>cd\).  On \(\mathcal G_\delta\),
\eqref{eq:core-degree-cost} and \(\delta\le1\) imply
\[
 d^{-j/2}|Z^{\Gamma_0}|
 \le d^{-j/2}(\delta\sqrt d)^D
 \le d^N\delta^j.
\]
Also, \(|\xi_e|\le1\) and
\(e^{B\sum_e|Z_e|^3/\sqrt d}\le e^{B\binom m2\delta^3d}\).  Hence the count
above gives
\[
 \sum_{\substack{\Gamma_0\in\mathfrak C_{N,s}(H)\\
                   2\omega(\Gamma_0)>cd}} I(\Gamma_0)
 \le C d^N e^{B\binom m2\delta^3d}
        \sum_{j>cd}(Q_m\delta)^j.
\]
Decrease \(\delta\), depending only on \(m,B\), so that
\(Q_m\delta<1\) and
\[
 c\log\frac1{Q_m\delta}>2B\binom m2\delta^3.
\]
The last display is then at most \(Cd^Ne^{-c'd}\) for some \(c'>0\).  Finally,
\[
 \frac{Cd^Ne^{-c'd}}{p^{h-C\eps_n}d^{-s}}
 =Cd^{N+s}\exp\left(-c'd+(h-C\eps_n)L\right)
 \le Cd^{N+K}e^{-c'd+KL}\longrightarrow0,
\]
because \(h,s\le K\), \(L=o(d)\), and \(\log d=o(d)\).  This proves
\eqref{eq:integrated-core-sum} and completes the proof.
\end{proof}

Lemma~\ref{claim:cost} below shows that any admissible product must pay
enough powers of \(d^{-1}\).  This is the bridge between the analytic Gram
expansion and the graph parameter \(\rho(H)\) in Proposition~\ref{prop:corr}.
Lemma~\ref{claim:cost} is purely combinatorial.  It says that centering forces
the monomial to cover each edge of \(H\), while symmetry forces even parity on
all auxiliary coordinates; those two requirements together force enough
generator cost.
Recall from Definition~\ref{def:admissible} that a generator product
\(\Gamma\) is \(H\)-admissible if
\(r_e(\Gamma)\ge1\) for every \(e\in E(H)\), and
\(r_e(\Gamma)\) is even for every
\(e\in E(K_V)\setminus E(H)\).
Recall also that the cost of a product
\(\Gamma=\{g_1,\ldots,g_r\}\), with multiplicities counted, is
\[
        \omega(\Gamma):=\sum_{a=1}^r\omega(g_a).
\]
Here a closed-walk generator of length \(\ell\) has cost
\(\ell/2-1\), while a one-edge correction generator \((e,s)\) has cost
\(1\) for \(s=1\) and \(s-1\) for \(s\ge2\).
\begin{lemma}[Covering cost]\label{claim:cost}
Let \(H\) be connected and have minimum degree at least \(2\).  Let
\(\Gamma\) be an \(H\)-admissible product of generators.  Then
\begin{equation}
        \omega(\Gamma)
        \ge
        \max\left\{
        \frac12,\,
        |V(H)|-1-\frac{|E(H)|}{2}
        \right\}.
        \label{eq:cover-cost}
\end{equation}
\end{lemma}

\begin{proof}
Write \(v=|V(H)|\) and \(h=|E(H)|\).
Because \(H\) has minimum degree at least \(2\), it is nonempty.
Since \(\Gamma\) is \(H\)-admissible, it covers every edge of \(H\), so
\(\Gamma\) contains at least one nonconstant generator.  A closed-walk
generator of length \(\ell\ge3\) has cost \(\ell/2-1\ge1/2\), whereas every
one-edge correction generator has cost at least \(1\).  Hence
\[
        \omega(\Gamma)\ge \frac12.
\]

\proofstep{Step 1: Decomposition into primitive pieces.}
Call a simple cycle of length at least \(3\), or two parallel copies of one
edge, a \emph{primitive piece}; the latter will be called a doubled edge.  For
multisets \(\mathcal C\) of cycles and \(\mathcal D\) of doubled edges, set
\[
 \operatorname{pcost}(\mathcal C,\mathcal D)
 :=\sum_{C\in\mathcal C}\left(\frac{|C|}{2}-1\right)
   +\frac{|\mathcal D|}{2}.
\]

The edge copies traversed by a closed-walk generator form a connected Eulerian
multigraph.  Repeatedly removing a shortest closed trail decomposes this
multigraph into primitive pieces: such a trail is either a simple cycle or a
doubled edge, and its removal preserves even degrees.  If a generator has
length \(\ell\) and its decomposition has \(a\) cycles and \(b\) doubled edges,
then its primitive cost is
\[
 \frac{\ell}{2}-a-\frac b2\le\frac{\ell}{2}-1=\omega(g).
\]
Indeed, either \(a\ge1\), or all pieces are doubled edges, in which case
\(\ell\ge3\) implies \(b\ge2\).  A one-edge correction \(Z_e^{2u}\) is replaced
by \(u\) doubled edges, whose total primitive cost \(u/2\) is at most
\(\omega(e,u)\): use \(1/2\le1\) for \(u=1\) and \(u/2\le u-1\) for
\(u\ge2\).  Taking the multiset union over all generators gives
primitive pieces with the same edge multiplicities as \(Z^\Gamma\) and total
primitive cost at most \(\omega(\Gamma)\).

\proofstep{Step 2: Removal of auxiliary edges.}
Every edge outside \(H\) has even multiplicity by \(H\)-admissibility.  Delete
any doubled edge outside \(H\).  If an auxiliary edge \(e=\{x,y\}\) remains,
it occurs in at least two cycles, counted with multiplicity.  Choose two such
cycles \(C_1,C_2\).  Removing \(e\) from each leaves two \(x\)-to-\(y\) paths
whose union is a connected Eulerian multigraph \(M\) of length
\[
        |E(M)|=|C_1|+|C_2|-2\ge4.
\]
Decomposing \(M\) as in Step~1 gives new primitive pieces of cost at most
\[
 \frac{|E(M)|}{2}-1
 =\left(\frac{|C_1|}{2}-1\right)
  +\left(\frac{|C_2|}{2}-1\right).
\]
This replacement removes two copies of \(e\), preserves every other edge
multiplicity, and does not increase the cost.  Repeating it, and deleting any
auxiliary doubled edges that arise, strictly decreases the total auxiliary-edge
multiplicity.  The procedure therefore terminates with primitive pieces
\((\mathcal C,\mathcal D)\) supported in \(H\), still covering \(H\), and
\[
        \operatorname{pcost}(\mathcal C,\mathcal D)\le\omega(\Gamma).
\]

\proofstep{Step 3: The cost estimate inside \(H\).}
List the cycles as \(C_1,\ldots,C_s\).  Let \(U_i\) be the simple graph formed
by the edges of \(C_1,\ldots,C_i\), and put
\[
 R_i:=\sum_{k=1}^i|C_k|,
 \qquad h_i:=|E(U_i)|,
 \qquad \beta(U_i):=|E(U_i)|-|V(U_i)|+c(U_i),
\]
with \(R_0=h_0=\beta(U_0)=0\).  The cyclomatic number \(\beta\) is monotone
under inclusion: adding an isolated vertex leaves it unchanged, while adding
an edge either leaves it unchanged or increases it by one.

For each \(i\),
\[
 1\le \beta(U_i)-\beta(U_{i-1})
 +\frac12\left[(R_i-h_i)-(R_{i-1}-h_{i-1})\right].
\]
To see this, if \(C_i\) contains a new edge \(e\), then \(e\) is not a bridge
of \(U_i\), and \(U_{i-1}\subseteq U_i-e\).  Hence
\(\beta(U_i)\ge\beta(U_{i-1})+1\), while the second term is nonnegative.  If
\(C_i\) contains no new edge, then \(U_i=U_{i-1}\) and the second term equals
\(|C_i|/2\ge3/2\).

Let \(R=R_s\) and \(h_{\mathcal C}=h_s\).  Summing the preceding inequality and
using \(U_s\subseteq H\) gives
\[
        s\le\beta(H)+\frac{R-h_{\mathcal C}}{2}.
\]
Every edge of \(H\) not covered by these cycles must be covered by a doubled
edge, so \(|\mathcal D|\ge h-h_{\mathcal C}\).  Since
\(\beta(H)=h-v+1\),
\[
\begin{aligned}
 \omega(\Gamma)
 &\ge \operatorname{pcost}(\mathcal C,\mathcal D)\\
 &=\frac R2-s+\frac{|\mathcal D|}{2}\\
 &\ge \frac R2-\beta(H)-\frac{R-h_{\mathcal C}}2
       +\frac{h-h_{\mathcal C}}2\\
 &=\frac h2-\beta(H)=v-1-\frac h2.
\end{aligned}
\]
Combining this with \(\omega(\Gamma)\ge1/2\) proves
\eqref{eq:cover-cost}.
\end{proof}

\subsection{Completion of the finite-graph correlation estimate}

\begin{proof}[Proof of the second assertion in Proposition~\ref{prop:corr}]
\leavevmode\par
Write \(h=|E(H)|\le K\) and
\(d\nu^\otimes(Z):=\prod_{e\in\mathcal E_m}d\nu_d(Z_e)\).  Fix
\(A=K+3\).  Let \(\delta=\delta_{m,C_m}\) be the cutoff supplied by
Lemma~\ref{lem:integrated-core} with \(B=C_m\).  By that lemma,
\(0<\delta\le\delta_m\), so all three assertions of
Lemma~\ref{lem:cluster} apply at this same cutoff.
We divide the argument into six steps.  The first two reduce the
correlation to expansion terms indexed by \(H\)-admissible products.  The next
two extract a bounded \(H\)-admissible subproduct from every such product, and
the final two group and sum the expansion according to these subproducts.

\proofstep{Step 1: Truncation and expansion.}
We first show that the complement of the convergence region is
negligible.
Since \(R_{m,d}\,d\nu^\otimes\) is the true joint law,
Equation~\eqref{eq:Phi-integral} and the third assertion of
Lemma~\ref{lem:cluster}, applied with \(W=1\) and \(F=E(H)\), give
\[
 \left|
 \int_{\mathcal G_\delta^c}
 \prod_{e\in E(H)}\xi_e(Z_e)R_{m,d}(Z)\,d\nu^\otimes(Z)
 \right|
 \le C_Kp^Ad^{-A}.
\]
As shown in the reduction to connected graphs, \(m\le h\le K\), and
therefore \(\rho(H)=\max\{1/2,m-1-h/2\}\le h\le K\).  Since
\(A=K+3\ge h,\rho(H)\), \(0<p<1\), and \(d\ge1\), the preceding error is at
most \(C_Kp^hd^{-\rho(H)}\).

It remains to analyze the integral on \(\mathcal G_\delta\).
On \(\mathcal G_\delta\), the first assertion of
Lemma~\ref{lem:cluster} gives
\[
 R_{m,d}(Z)=\sum_\Gamma
 c_\Gamma d^{-\omega(\Gamma)}Z^\Gamma.
\]
The absolute-sum bound \eqref{eq:absolute-subproduct-bound} in
Lemma~\ref{lem:cluster}, with \(\Gamma_0=\varnothing\), and
\eqref{eq:weighted-uniform-full} in Lemma~\ref{lem:integrated-core}, with
\(r=0\), imply
\[
 \int_{\mathcal G_\delta}
 \sum_\Gamma |c_\Gamma|d^{-\omega(\Gamma)}|Z^\Gamma|
 \,d\nu^\otimes(Z)
 \le C_K.
\]
Since \(|\xi_e|\le1\), Tonelli's theorem therefore permits termwise
integration on \(\mathcal G_\delta\).

\proofstep{Step 2: Eliminating non-admissible products.}
Suppose first that \(r_f(\Gamma)\) is odd for some
\(f\in\mathcal E_m\setminus E(H)\).  Then
\[
 \int_{\mathcal G_\delta}
 \prod_{e\in E(H)}\xi_e(Z_e)Z^\Gamma\,d\nu^\otimes(Z)=0,
\]
because both \(\nu_d\) and the truncation interval
\([-\delta\sqrt d,\delta\sqrt d]\) are symmetric in \(Z_f\).

We next consider products that omit an edge of \(H\).  Without
truncation, their integrals vanish because \(\int\xi_e\,d\nu_d=0\); after
truncation, the resulting error is exponentially small.  Fix \(e\in E(H)\), and set
\[
 \zeta_e:=\int_{|z|\le\delta\sqrt d}\xi_e(z)\,d\nu_d(z).
\]
Since \(\int\xi_e\,d\nu_d=0\) and \(|\xi_e|\le1\),
Lemma~\ref{lem:tail}, specifically \eqref{eq:basic-gaussian-tail}, gives
\[
 |\zeta_e|
 \le \Pp\{|Z|>\delta\sqrt d\}
 \le C_Ke^{-c_Kd}.
\]
Let \(\mathcal G_\delta^{(-e)}\) be the truncation event for all coordinates
other than \(e\), and write \(d\nu_{-e}^\otimes\) for their product law.
For \(Z_{-e}\in\mathcal G_\delta^{(-e)}\), set
\(\widetilde Z_e=0\) and \(\widetilde Z_f=Z_f\) for \(f\ne e\).
Evaluating the full series at \(\widetilde Z\) removes exactly the
terms for which \(r_e(\Gamma)>0\).
The absolute-sum bound \eqref{eq:absolute-subproduct-bound} in
Lemma~\ref{lem:cluster}, again with the empty subproduct, yields
\[
\begin{aligned}
 \sum_{\Gamma:r_e(\Gamma)=0}
 |c_\Gamma|d^{-\omega(\Gamma)}|Z^\Gamma|
 &=
 \sum_\Gamma |c_\Gamma|d^{-\omega(\Gamma)}
 |\widetilde Z^\Gamma|\\
 &\le C_K\exp\left(
 C_m\sum_{f\ne e}\frac{|Z_f|^3}{\sqrt d}\right).
\end{aligned}
\]
For every term with \(r_e(\Gamma)=0\), the monomial \(Z^\Gamma\) is
independent of \(Z_e\); hence the integral in the \(e\)-coordinate is
\(\zeta_e\).  Using the product structure of \(\mathcal G_\delta\), followed by
\eqref{eq:weighted-uniform-centered} and
\eqref{eq:weighted-uniform-full} in Lemma~\ref{lem:integrated-core}, we obtain
\[
\begin{aligned}
 S_e
 &:=
 \sum_{\Gamma:r_e(\Gamma)=0}
 |c_\Gamma|d^{-\omega(\Gamma)}
 \left|
 \int_{\mathcal G_\delta}
 \prod_{f\in E(H)}\xi_f(Z_f)Z^\Gamma\,d\nu^\otimes(Z)
 \right|\\
 &\le
 |\zeta_e|
 \int_{\mathcal G_\delta^{(-e)}}
 \prod_{f\in E(H)\setminus\{e\}}|\xi_f(Z_f)|
 \sum_{\Gamma:r_e(\Gamma)=0}
 |c_\Gamma|d^{-\omega(\Gamma)}|Z^\Gamma|
 \,d\nu_{-e}^\otimes(Z_{-e})\\
 &\le C_Ke^{-c_Kd}p^{h-1-C_K\eps_n}.
\end{aligned}
\]
Moreover,
\[
 \frac{e^{-c_Kd}p^{h-1-C_K\eps_n}}{p^Ad^{-A}}
 =
 \exp\left(-c_Kd+O_K(L+\eps_nL+\log d)\right)
 \longrightarrow0,
\]
because \(L=o(d)\), \(\eps_nL=o(d)\), and \(\log d=o(d)\).  Summing \(S_e\)
over the at most \(K\) edges of \(H\), and combining this bound with the exact
parity cancellation, gives
\[
\begin{aligned}
 &\int_{\mathcal G_\delta}
 \prod_{e\in E(H)}\xi_e(Z_e)R_{m,d}(Z)\,d\nu^\otimes(Z)\\
 &\qquad=
 \sum_{\Gamma:\,\Gamma\text{ is }H\text{-admissible}}
 c_\Gamma d^{-\omega(\Gamma)}
 \int_{\mathcal G_\delta}
 \prod_{e\in E(H)}\xi_e(Z_e)Z^\Gamma\,d\nu^\otimes(Z)
 +O_K(p^Ad^{-A}).
\end{aligned}
\]

\proofstep{Step 3: Covering \(H\) and identifying the parity defect.}

An \(H\)-admissible product may contain arbitrarily many generators, whereas
Lemma~\ref{lem:integrated-core} applies to subproducts with a fixed bound on
their number of factors.  Recall that \(h=|E(H)|\), set
\(N_K:=\binom K2\), and recall that \(\mathfrak C_{N,s}(H)\) is the set of
\(H\)-admissible generator products \(\Gamma_0\) such that

\[
        |\Gamma_0|\le N
        \qquad\text{and}\qquad
        \omega(\Gamma_0)\ge s.
\]

We will prove that every \(H\)-admissible product \(\Gamma\) contains a
subproduct in \(\mathfrak C_{N_K,\rho(H)}(H)\).  The construction has two
parts: we first select enough factors to cover \(H\), and then add a bounded
number of factors to make every auxiliary-edge exponent even.

Set \(\mathcal E_{\rm aux}:=\mathcal E_m\setminus E(H)\).  For a generator
subproduct \(\Lambda\), let

\[
 \mathcal O(\Lambda)
 :=\{f\in\mathcal E_{\rm aux}:r_f(\Lambda)\text{ is odd}\}.
\]
Thus \(\mathcal O(\Lambda)\) records exactly the auxiliary-edge parity defect
of \(\Lambda\); it is empty precisely when every auxiliary-edge exponent is
even.

Fix an \(H\)-admissible \(\Gamma\), and temporarily distinguish its generator
factors even when some of them are identical.

For each \(e\in E(H)\), choose one factor whose monomial contains \(Z_e\).
Such a factor exists because
\(r_e(\Gamma)\ge1\).  Let \(\Gamma_{\rm cov}\) consist of all factors selected
in this way, counting a factor only once if it was selected for several edges.
Then

\[
        \Gamma_{\rm cov}\text{ covers }H,
        \qquad |\Gamma_{\rm cov}|\le h.
\]

Let \(\Gamma_{\rm rem}\) consist of the factors of \(\Gamma\) not selected for
\(\Gamma_{\rm cov}\).  Thus
\(\Gamma=\Gamma_{\rm cov}\uplus\Gamma_{\rm rem}\).  Since \(\Gamma\) is
\(H\)-admissible, \(r_f(\Gamma)\) is even for every
\(f\in\mathcal E_{\rm aux}\).  Consequently,

\[
 r_f(\Gamma_{\rm rem})
 =r_f(\Gamma)-r_f(\Gamma_{\rm cov})
 \equiv r_f(\Gamma_{\rm cov})\pmod 2,
 \qquad f\in\mathcal E_{\rm aux}.
\]

Equivalently,
\(\mathcal O(\Gamma_{\rm rem})=\mathcal O(\Gamma_{\rm cov})\).  Thus the
whole remaining product corrects the auxiliary-edge parities introduced by
\(\Gamma_{\rm cov}\).  The next step shows that only a bounded part of
\(\Gamma_{\rm rem}\) is needed.

\proofstep{Step 4: Correcting the parity with few factors.}

We now reproduce the parity defect of \(\Gamma_{\rm cov}\) using as few
factors from \(\Gamma_{\rm rem}\) as possible.  Among all subproducts
\(\Lambda\subseteq\Gamma_{\rm rem}\) satisfying

\[
        \mathcal O(\Lambda)=\mathcal O(\Gamma_{\rm cov}),
\]

choose one with the fewest factors and call it \(\Gamma_{\rm par}\).  This
choice is possible because \(\Gamma_{\rm rem}\) itself satisfies the equality.
Write \(t:=|\Gamma_{\rm par}|\) and
\(N_{\rm aux}:=|\mathcal E_{\rm aux}|\).  We claim that
\(t\le N_{\rm aux}\).

To prove the claim, label the factors of \(\Gamma_{\rm par}\) as
\(g_1,\ldots,g_t\).  For \(I\subseteq[t]\), let \(\Gamma_I\) be the subproduct
formed by the factors \(g_i\) with \(i\in I\).  Minimality first implies

\[
 \varnothing\ne D\subseteq[t]
 \qquad\Longrightarrow\qquad
 \mathcal O(\Gamma_D)\ne\varnothing.
\]

Indeed, if \(\mathcal O(\Gamma_D)=\varnothing\) for some nonempty
\(D\subseteq[t]\), then, for every auxiliary edge \(f\),

\[
 r_f(\Gamma_{[t]\setminus D})
 \equiv r_f(\Gamma_{\rm par})-r_f(\Gamma_D)
 \equiv r_f(\Gamma_{\rm par})\pmod2.
\]

The complementary subproduct \(\Gamma_{[t]\setminus D}\) would then have the
same odd-edge set as \(\Gamma_{\rm par}\) and fewer factors, contradicting
minimality.

We can now count parity patterns.  If distinct sets \(I,J\subseteq[t]\)
satisfied
\(\mathcal O(\Gamma_I)=\mathcal O(\Gamma_J)\), then
\(D:=I\triangle J\) would be nonempty and, for every auxiliary edge \(f\),

\[
 r_f(\Gamma_D)
 \equiv r_f(\Gamma_I)+r_f(\Gamma_J)
 \equiv0\pmod2.
\]

This would give \(\mathcal O(\Gamma_D)=\varnothing\), contrary to the preceding
implication.  Thus the map
\(I\mapsto\mathcal O(\Gamma_I)\) is injective.  Since an odd-edge set is a
subset of the \(N_{\rm aux}\) auxiliary edges, there are only
\(2^{N_{\rm aux}}\) possible such sets.  Therefore

\[
        2^t\le2^{N_{\rm aux}},
        \qquad\text{so}\qquad
        |\Gamma_{\rm par}|=t\le N_{\rm aux}.
\]

Now set \(\Gamma_0:=\Gamma_{\rm cov}\uplus\Gamma_{\rm par}\).  The first
subproduct covers \(H\), while the second corrects its auxiliary-edge parities.
More precisely, for every \(f\in\mathcal E_{\rm aux}\),

\[
 r_f(\Gamma_0)
 \equiv r_f(\Gamma_{\rm cov})+r_f(\Gamma_{\rm par})
 \equiv0\pmod2.
\]

Hence \(\Gamma_0\) is \(H\)-admissible.  Its number of factors satisfies

\[
 |\Gamma_0|
 \le h+N_{\rm aux}
 =|E(H)|+|\mathcal E_m\setminus E(H)|
 =\binom m2\le N_K.
\]

This uniform bound, independent of \(\Gamma\), allows us to apply the
subproduct bound in Lemma~\ref{lem:cluster} and the integrated summation in
Lemma~\ref{lem:integrated-core} with the fixed choice \(N=N_K\) in
Steps~5--6.  Finally, Lemma~\ref{claim:cost}, specifically
\eqref{eq:cover-cost}, gives \(\omega(\Gamma_0)\ge\rho(H)\).  Therefore
\(\Gamma_0\in\mathfrak C_{N_K,\rho(H)}(H)\), as required.

\proofstep{Step 5: Overcounting by bounded subproducts.}
Steps~3 and~4 show that every \(H\)-admissible \(\Gamma\) contains at least one
\(\Gamma_0\in\mathfrak C_{N_K,\rho(H)}(H)\).  Thus, for every generator product
\(\Gamma\),
\[
 \1_{\{\Gamma\text{ is }H\text{-admissible}\}}
 \le
 \sum_{\Gamma_0\in\mathfrak C_{N_K,\rho(H)}(H)}
 \1_{\{\Gamma\supseteq\Gamma_0\}}.
\]
No choice of a distinguished subproduct is needed; if \(\Gamma\) contains more
than one such \(\Gamma_0\), the right-hand side simply overcounts it.  Therefore,
for \(Z\in\mathcal G_\delta\), Lemma~\ref{lem:cluster}, specifically
\eqref{eq:absolute-subproduct-bound}, gives
\[
\begin{aligned}
 &\sum_{\Gamma:\,\Gamma\text{ is }H\text{-admissible}}
 |c_\Gamma|d^{-\omega(\Gamma)}|Z^\Gamma|\\
 &\quad\le
 \sum_{\Gamma_0\in\mathfrak C_{N_K,\rho(H)}(H)}
 \sum_{\Gamma\supseteq\Gamma_0}
 |c_\Gamma|d^{-\omega(\Gamma)}|Z^\Gamma|\\
 &\quad\le
 C_K\sum_{\Gamma_0\in\mathfrak C_{N_K,\rho(H)}(H)}
 d^{-\omega(\Gamma_0)}|Z^{\Gamma_0}|
 \exp\left(C_m\sum_{f\in\mathcal E_m}
 \frac{|Z_f|^3}{\sqrt d}\right).
\end{aligned}
\]
Here the constant in Lemma~\ref{lem:cluster} is absorbed into \(C_K\), because
\(m\le K\) and \(N_K=\binom K2\).

\proofstep{Step 6: Summation over the bounded subproducts.}
Define
\[
 S_H:=
 \sum_{\Gamma:\,\Gamma\text{ is }H\text{-admissible}}
 c_\Gamma d^{-\omega(\Gamma)}
 \int_{\mathcal G_\delta}
 \prod_{e\in E(H)}\xi_e(Z_e)Z^\Gamma\,d\nu^\otimes(Z).
\]
The absolute-integrability bound from Step~1 permits the triangle inequality
and termwise integration.  Applying the pointwise bound from Step~5 and then
Tonelli's theorem gives
\[
\begin{aligned}
 |S_H|
 &\le C_K
 \sum_{\Gamma_0\in\mathfrak C_{N_K,\rho(H)}(H)}
 d^{-\omega(\Gamma_0)}
 \int_{\mathcal G_\delta}
 \prod_{e\in E(H)}|\xi_e(Z_e)|\,|Z^{\Gamma_0}|
 \exp\left(C_m\sum_{f\in\mathcal E_m}
 \frac{|Z_f|^3}{\sqrt d}\right)
 \,d\nu^\otimes(Z).
\end{aligned}
\]
Finally, \eqref{eq:integrated-core-sum} in
Lemma~\ref{lem:integrated-core}, with \(N=N_K\) and \(s=\rho(H)\), bounds the
last expression by
\[
 C_Kp^{h-C_K\eps_n}L^{C_K}d^{-\rho(H)}.
\]

Combining the error estimates from Steps~1 and~2 with the bound on
\(S_H\) and using
\[
 p^Ad^{-A}
 \le p^hd^{-\rho(H)}
 \le p^{h-C_K\eps_n}L^{C_K}d^{-\rho(H)}
\]
for all sufficiently large \(n\), proves
\eqref{eq:correlation-bound} after enlarging the constants.  This completes
the proof of Proposition~\ref{prop:corr}.
\end{proof}

\section{Proof of Theorem~\ref{thm:main}}\label{sec:proof_semicircle}

We now combine the finite-graph correlation estimate with the trace-moment method.
The proof separates the usual Wigner tree contribution from all terms that
contain a nonempty centered graph, and then proves concentration of the moments.
Under the assumptions of Theorem~\ref{thm:main},
\(L/d\le 1/(c_0np)\to0\), so Proposition~\ref{prop:corr} applies and
\(\eps_n=\sqrt{L/d}=o(1)\).

\subsection{Reduction of repeated edge powers}\label{subsec:power-reduction}
The moment expansion contains powers of the same centered edge variable whenever
a closed walk traverses an edge more than once.
Lemma~\ref{lem:power-reduction} reduces every such power to a centered part
plus a scalar part.

\begin{lemma}[Repeated centered-edge powers]\label{lem:power-reduction}
Let \(U,V\) be independent uniform points on \(S^{d-1}\), and set
\(Y=\xi(U,V)=\1_{\{\ip{U}{V}\ge\tau\}}-p\).
Then \(Y=q\) with probability \(p\), and \(Y=-p\) with probability \(q=1-p\).
For every integer \(r\ge2\), there are constants \(\alpha_r,\beta_r\) such that
\begin{equation}
        Y^r=\alpha_rY+\beta_r.
        \label{eq:power-reduction}
\end{equation}
Moreover,
\begin{equation}
        |\alpha_r|\le 1,
        \qquad
        |\beta_r|\le p,
        \label{eq:alpha-beta-bound}
\end{equation}
and
\begin{equation}
        \beta_2=pq.
        \label{eq:beta2}
\end{equation}
\end{lemma}

\begin{proof}
Set
\[
        \alpha_r=q^r-(-p)^r,
        \qquad
        \beta_r=pq^r+q(-p)^r.
\]
Since \(Y\) takes only the values \(q\) and \(-p\), direct substitution, using
\(q+p=1\), verifies \eqref{eq:power-reduction} with the displayed coefficients.
Taking expectations and using \(\E Y=0\) gives \(\beta_r=\E Y^r\).  Moreover,
\(|\alpha_r|\le q^r+p^r\le1\), and, for \(r\ge2\),
\[
        |\beta_r|\le pq^r+qp^r
        =pq(q^{r-1}+p^{r-1})\le pq\le p.
\]
The case \(r=2\) gives \(\beta_2=pq(q+p)=pq\).
\end{proof}

\subsection{Expected trace moments}
We now apply the trace-moment method.  The role of
Lemma~\ref{lem:power-reduction} is to separate each walk contribution into a
purely scalar part, which will produce the Catalan moments, and a centered graph
part, which will be killed or bounded by Proposition~\ref{prop:corr}.

Fix \(k\ge1\).  Since \((M_n)_{ii}=0\),
\begin{equation}
        \E\frac1n\Tr M_n^k
        =\frac{1}{n(npq)^{k/2}}
        \sum_{i_1,\ldots,i_k}
        \E\prod_{\ell=1}^k\xi_{i_\ell i_{\ell+1}},
        \label{eq:moment-expansion}
\end{equation}
where \(i_{k+1}=i_1\), and summands with \(i_\ell=i_{\ell+1}\) vanish.

A cyclic word \(w=(i_1,\ldots,i_k)\) defines a connected multigraph \(G_w\).
Let \(v(w)\) be the number of distinct vertices visited by the word, let
\(e(w)\) be the number of distinct unoriented edges used by the word, and let
\(m_f(w)\) be the multiplicity of edge \(f\).  For a fixed unlabeled word
pattern, the number of labelings is \(O_k(n^{v(w)})\).  When \(w\) is fixed, we
abbreviate \(v=v(w)\), \(e=e(w)\), and \(m_f=m_f(w)\).

For a fixed labeled cyclic word \(w\), group the product in \eqref{eq:moment-expansion}
according to the distinct unoriented edges of \(G_w\):
\[
        \prod_{\ell=1}^k \xi_{i_\ell i_{\ell+1}}
        =
        \prod_{f\in E(G_w)} \xi_f^{m_f(w)}.
\]
Edges with multiplicity \(1\) already appear as a single centered factor
\(\xi_f\), so they are not reduced and must remain centered.  For every edge
\(f\) with multiplicity \(m_f(w)\ge2\), apply
Lemma~\ref{lem:power-reduction}, specifically
\eqref{eq:power-reduction}:
\[
        \xi_f^{m_f}
        =
        \alpha_{m_f}\xi_f+\beta_{m_f}.
\]
Thus, in each term of the resulting expansion, every repeated edge makes one
of two choices: \(\alpha_{m_f}\xi_f\) or \(\beta_{m_f}\).
The set of repeated edges is
\(R(w):=\{f\in E(G_w):m_f(w)\ge2\}\).
For each expansion term, write
\(R(w)=S\sqcup C\),
where \(S\) is the set of repeated edges for which we choose the centered term
\(\alpha_{m_f}\xi_f\), and \(C\) is the set of repeated edges for which we
choose the constant term \(\beta_{m_f}\).
Define the simple graph \(H=H(w,S)\) of centered edges by
\(E(H)=\{f\in E(G_w):m_f(w)=1\}\cup S\).
Its vertex set is the set of endpoints of these edges.  Thus every
edge of multiplicity one is automatically included in \(H\), while a repeated edge is
included in \(H\) only if its centered part is chosen.

With this notation, each expansion term has the form
\begin{equation}
        \left(\prod_{f\in C}\beta_{m_f}\right)
        \left(\prod_{f\in S}\alpha_{m_f}\right)
        \E\prod_{f\in E(H)}\xi_f.
        \label{eq:expanded-term}
\end{equation}
If \(H\) has a leaf, then the last expectation is zero by
Proposition~\ref{prop:corr}.  We therefore only need to consider the two
remaining cases: \(H=\varnothing\), and \(H\ne\varnothing\) with every
component of \(H\) having minimum degree at least \(2\).
\subsection{Pure constant terms and Catalan moments}
We first identify the terms in which every repeated edge has been replaced by
its scalar part.  These are the only terms capable of matching the usual Wigner
tree contribution.

If \(H=\varnothing\), then every edge in the word has multiplicity at least \(2\)
and has been replaced by a constant.  Using
Lemma~\ref{lem:power-reduction}, specifically \eqref{eq:alpha-beta-bound}, the
contribution of a fixed word shape is bounded by
\begin{equation}
        C_k\frac{n^v p^e}{n(npq)^{k/2}}
        \le C_kn^{v-1-k/2}p^{e-k/2}.
        \label{eq:constant-bound}
\end{equation}
Since the support is connected, \(v\le e+1\).  Since every edge is repeated,
\(e\le k/2\).  If \(e<k/2\), then \eqref{eq:constant-bound} is at most
\(C_k(np)^{e-k/2}=o(1)\).
If \(e=k/2\) but \(v<e+1\), the same bound is
\(O_k(n^{-1})=o(1)\).
The only possible leading case is
\[
        e=k/2,
        \qquad
        v=e+1.
\]
Then the support graph is a tree and every edge is traversed exactly twice.
These are the usual Wigner tree words.  Such words exist only for even \(k\).
After relabeling vertices in order of first appearance, such a word
is the contour traversal of a rooted plane tree.  This is the standard bijection
that produces Catalan numbers in the moment proof of the semicircle law.
Indeed, if \(k=2r\), then the leading case has
\(e=r\) and \(v=r+1\).
For each Wigner tree word, there are
\(n(n-1)\cdots(n-r)=(1+o(1))n^{r+1}\) labelings.  Since every edge is traversed
exactly twice, \eqref{eq:beta2} in Lemma~\ref{lem:power-reduction} shows that
each edge contributes \(\beta_2=pq\).  Hence the contribution of one such word is
\[
        \frac{(1+o(1))n^{r+1}(pq)^r}{n(npq)^r}
        =
        1+o(1).
\]
The number of such word patterns is the Catalan number \(\Cat_r\).  Therefore the
total leading pure constant contribution is
\(\Cat_r(1+o(1))=\Cat_r+o(1)\).
If \(k\) is odd, no integer \(e\) can satisfy \(e=k/2\), so all pure constant
terms are \(o(1)\).  Thus the contribution of all pure constant terms is
\begin{equation}
        \begin{cases}
        \Cat_{k/2}+o(1),& k\text{ even},\\[1mm]
        o(1),& k\text{ odd}.
        \end{cases}
        \label{eq:wigner-contribution}
\end{equation}
\subsection{Terms with a nonempty centered graph}
It remains to show that every term with a nonempty centered graph is
negligible.  The argument combines the correlation bound with two elementary
counting inequalities for the support of the word.

Now suppose \(H\ne\varnothing\) and no component of \(H\) has a leaf.  Let
\[
        h=|E(H)|,
        \qquad
        u=|V(H)|,
        \qquad
        a_H=\#\Comp(H),
        \qquad
        c=|C|.
\]
Contract each of the \(a_H\) components of \(H\) to one vertex and
leave the \(v-u\) vertices outside \(H\) unchanged.  Because the full word support is
connected, the resulting graph is connected; all of its nonloop edges come
from the \(c\) constant edges.  A connected graph on
\(a_H+v-u\) vertices has at least \(a_H+v-u-1\) edges, and therefore
\[
        c\ge a_H+v-u-1.
\]
Equivalently,
\begin{equation}
        v\le u+c-a_H+1.
        \label{eq:v-bound}
\end{equation}
Every centered edge appears at least once in the word, while every constant edge
appears at least twice.  Therefore
\begin{equation}
        k\ge h+2c.
        \label{eq:k-bound}
\end{equation}
By Proposition~\ref{prop:corr} and Lemma~\ref{lem:power-reduction},
specifically \eqref{eq:alpha-beta-bound}, the absolute contribution of this
fixed shape is bounded by
\begin{equation*}
        C_k\frac{n^v p^c p^{h-C_k\eps_n}L^{C_k}d^{-\rho(H)}}{n(npq)^{k/2}}.
\end{equation*}
Since \(q\ge 1/2\), we have \(q^{-k/2}=O_k(1)\), so the preceding expression is bounded by
\begin{equation*}
        C_k p^{-C_k\eps_n}L^{C_k}
        n^{v-1-k/2}p^{c+h-k/2}d^{-\rho(H)}.
\end{equation*}
Define
\[
        \Delta_w:=\frac{k}{2}-\left(\frac h2+c\right).
\]
By \eqref{eq:k-bound}, \(k\ge h+2c\), and hence \(\Delta_w\ge0\).  Equivalently,
\[
        \frac{k}{2}=\frac h2+c+\Delta_w .
\]
Using \eqref{eq:v-bound}, we obtain
\[
\begin{aligned}
        v-1-\frac{k}{2}
        &\le
        u+c-a_H-\frac{k}{2} \\
        &=
        u+c-a_H-\left(\frac h2+c+\Delta_w\right) \\
        &=
        u-a_H-\frac h2-\Delta_w.
\end{aligned}
\]
Moreover,
\[
        c+h-\frac{k}{2}
        =
        c+h-\left(\frac h2+c+\Delta_w\right)
        =
        \frac h2-\Delta_w.
\]
Therefore,
\[
\begin{aligned}
        n^{v-1-k/2}p^{c+h-k/2}
        &\le
        n^{u-a_H-h/2-\Delta_w}p^{h/2-\Delta_w} \\
        &=
        n^{u-a_H-h/2}p^{h/2}(np)^{-\Delta_w} \\
        &\le
        n^{u-a_H-h/2}p^{h/2},
\end{aligned}
\]
where the last inequality holds for all sufficiently large \(n\), since
\(np\to\infty\) and \(\Delta_w\ge0\).  Hence
\begin{equation}
        C_k
        \frac{n^v p^c p^{h-C_k\eps_n}L^{C_k}d^{-\rho(H)}}
        {n(npq)^{k/2}}
        \le
        C_kp^{-C_k\eps_n}L^{C_k}
        n^\theta p^{h/2}d^{-\rho(H)},
        \label{eq:centered-main-bound}
\end{equation}
where \(\theta=u-a_H-h/2\).
We now show that the right-hand side of
\eqref{eq:centered-main-bound} is \(o(1)\).  We split according to the sign of
\(\theta\).

\emph{Case 1: \(\theta>0\).} From the definition of \(\rho(H)\),
we have \(\rho(H)\ge\theta\).
Since \(d\ge c_0npL\),
\[
\begin{aligned}
        n^\theta d^{-\rho(H)}
        \le n^\theta d^{-\theta}
        =\left(\frac nd\right)^\theta
        \le c_0^{-\theta}(pL)^{-\theta}.
\end{aligned}
\]
Since \(k\) is fixed, the factor \(c_0^{-\theta}\) may be absorbed into
\(C_{k,c_0}\).  Thus \eqref{eq:centered-main-bound} is bounded by
\begin{equation}
        C_{k,c_0}p^{h/2-\theta-C_k\eps_n}L^{C_k-\theta}.
        \label{eq:theta-positive-bound}
\end{equation}
Moreover,
\begin{equation}
\begin{aligned}
        \frac h2-\theta
        =h-u+a_H
        = \sum_{D\in\Comp(H)}
        \bigl(|E(D)|-|V(D)|+1\bigr).
\end{aligned}
        \label{eq:cycle-rank}
\end{equation}
Every component \(D\) of \(H\) has minimum degree at least \(2\).  Thus
\[
        2|E(D)|
        =
        \sum_{v\in V(D)}\deg_D(v)
        \ge
        2|V(D)|.
\]
Thus
\begin{equation}
        \frac h2-\theta\ge a_H\ge1.
        \label{eq:cycle-rank-positive}
\end{equation}
Since \(\eps_n\to0\), for all sufficiently large \(n\),
\(h/2-\theta-C_k\eps_n\ge1/2\).
Therefore
\[
        \eqref{eq:theta-positive-bound}
        \le
        C_{k,c_0}p^{1/2}L^{C_k-\theta}
        =
        o(1),
\]
because \(p^\gamma(\log(1/p))^C\to0\) for every fixed
\(\gamma>0\) and \(C\in\R\).

\emph{Case 2: \(\theta \le 0\)}.  Then \(n^\theta\le1\).
Since \(H\ne\varnothing\), we have
\(\rho(H)\ge1/2\).
Furthermore, \(H\) is a nonempty simple graph and every component has minimum
degree at least \(2\).  Hence \(h=|E(H)|\ge3\).  Thus
\eqref{eq:centered-main-bound} is at most
\[
        C_kp^{h/2-C_k\eps_n}L^{C_k}d^{-1/2}.
\]
Since \(\eps_n\to0\), for all sufficiently large \(n\),
\(h/2-C_k\eps_n\ge1\).
Therefore
\[
        C_kp^{h/2-C_k\eps_n}L^{C_k}d^{-1/2}
        \le
        C_kpL^{C_k}d^{-1/2}
        =
        o(1).
\]
Here \(pL^{C_k}=e^{-L}L^{C_k}\to0\), while
\(d\ge c_0npL\to\infty\).

Since \(k\) is fixed, there are only finitely many unlabeled word patterns
and finitely many choices in the repeated-edge expansion.  Hence the sum of
all terms with a nonempty centered graph is \(o(1)\).  Combining this with
\eqref{eq:wigner-contribution}, we obtain
\begin{equation}
        \E\frac1n\Tr M_n^k
        \longrightarrow
        \begin{cases}
        0,& k\text{ odd},\\[1mm]
        \Cat_{k/2},& k\text{ even}.
        \end{cases}
        \label{eq:expected-moment-limit}
\end{equation}
\subsection{Variance estimate}
The same word estimates also give concentration of the trace moments.  We
first obtain a uniform bound for one word pattern and then sum over the finitely
many possible overlap patterns.

Let \(m_k(n)=n^{-1}\Tr M_n^k\).  For a labeled cyclic word
\(w=(i_1,\ldots,i_k)\), set \(i_{k+1}=i_1\), assume
\(i_\ell\ne i_{\ell+1}\), and define
\[
        X(w)=\prod_{\ell=1}^k\xi_{i_\ell i_{\ell+1}},
        \qquad
        V(w)=\{i_1,\ldots,i_k\}.
\]
Then
\[
        m_k(n)
        =
        \frac{1}{n(npq)^{k/2}}
        \sum_w X(w),
\]
and therefore
\[
        \Var(m_k(n))
        =
        \frac{1}{n^2(npq)^k}
        \sum_{w_1,w_2}
        \operatorname{Cov}\bigl(X(w_1),X(w_2)\bigr),
\]
where both sums range over labeled cyclic words of length \(k\).  If
\(V(w_1)\cap V(w_2)=\varnothing\), then \(X(w_1)\) and \(X(w_2)\) depend on
disjoint families of independent sphere points, so their covariance is zero.
It therefore remains to consider pairs of words sharing at least one vertex.
By the triangle inequality,
\begin{equation}
        \Var(m_k(n))\le V_1+V_2,
        \label{eq:variance-split}
\end{equation}
where
\[
        V_1
        :=
        \frac{1}{n^2(npq)^k}
        \sum_{\substack{w_1,w_2\\
        V(w_1)\cap V(w_2)\ne\varnothing}}
        \left|\E[X(w_1)X(w_2)]\right|
\]
and
\[
        V_2
        :=
        \frac{1}{n^2(npq)^k}
        \sum_{\substack{w_1,w_2\\
        V(w_1)\cap V(w_2)\ne\varnothing}}
        \left|\E X(w_1)\right|
        \left|\E X(w_2)\right|.
\]
We use the following consequence of the expected-moment estimates above.  Let
\(P\) be a fixed
unlabeled cyclic word pattern of length \(\ell\in\{k,2k\}\), and let \(v(P)\)
be the number of vertices of \(P\).  If \(w\) is any injective labeling of
\(P\), then
\begin{equation}
        \frac{n^{v(P)}}{n(npq)^{\ell/2}}
        \left|\E X(w)\right|
        \le C_{k,c_0}.
        \label{eq:word-pattern-bound}
\end{equation}
By exchangeability, \(\E X(w)\) is the same for every injective
labeling of \(P\), and \(P\) has \(\Theta_P(n^{v(P)})\) such labelings.  Thus
the normalization in \eqref{eq:word-pattern-bound} is, up to a constant
depending only on \(P\), the normalization used above for the total contribution
of that pattern.  To prove the bound, repeat the expansion used for the expectation with the single fixed
pattern \(P\).  If the centered graph produced by the expansion has a leaf, the
term vanishes.  If it is empty, the estimate is exactly the pure-constant bound
\eqref{eq:constant-bound}.  If it is nonempty and has no leaf, the derivation of
\eqref{eq:centered-main-bound} and the two cases following it give a bound
\(O_{k,c_0}(1)\) after multiplication by \(n^{v(P)}/(n(npq)^{\ell/2})\).  Since
\(\ell\le2k\), there are only finitely many choices in the repeated-edge
expansion, and \eqref{eq:word-pattern-bound} follows.

\paragraph{\texorpdfstring{Estimating \(V_1\)}{Estimating V1}}

Fix an overlap pattern for two cyclic words \(P_1,P_2\) of length \(k\).  Let
\(v_j=|V(P_j)|\) and \(r=|V(P_1)\cap V(P_2)|\ge1\).  Their union has
\(v_{12}=v_1+v_2-r\) vertices, so the pattern has
\(O_k(n^{v_{12}})\) labelings.  For one such labeling, cyclically rotate
\(w_1,w_2\) to begin at a common vertex and concatenate them.  The resulting
cyclic word \(w_{12}\) has length \(2k\), has \(v_{12}\) vertices, and satisfies
\(X(w_{12})=X(w_1)X(w_2)\).  Therefore \eqref{eq:word-pattern-bound} gives
\[
\begin{aligned}
        \frac{C_{k,c_0} n^{v_{12}}}{n^2(npq)^k}
        \left|\E[X(w_1)X(w_2)]\right|
        &=\frac{C_{k,c_0}}{n}
        \left(
        \frac{n^{v_{12}}}{n(npq)^k}
        \left|\E X(w_{12})\right|
        \right) \\
        &\le \frac{C_{k,c_0}}n=o(1).
\end{aligned}
\]
Summing over the finitely many overlap patterns gives
\begin{equation}
        V_1=o(1).
        \label{eq:V1-small}
\end{equation}

\paragraph{\texorpdfstring{Estimating \(V_2\)}{Estimating V2}}

For the same overlap pattern, apply \eqref{eq:word-pattern-bound} separately to
\(w_1\) and \(w_2\).  Since the number of labelings is
\(O_k(n^{v_1+v_2-r})\), its contribution to \(V_2\) is at most
\[
\begin{aligned}
        &C_{k,c_0}\frac{
        n^{v_1+v_2-r}
        \left|\E X(w_1)\right|
        \left|\E X(w_2)\right|
        }{n^2(npq)^k} \\
        &\qquad=
        C_{k,c_0}n^{-r}
        \prod_{j=1}^2\left(
        \frac{n^{v_j}}
        {n(npq)^{k/2}}
        \left|\E X(w_j)\right|
        \right) \\
        &\qquad\le C_{k,c_0}n^{-r}
        \le C_{k,c_0}n^{-1}=o(1).
\end{aligned}
\]
Here the last inequality uses \(r\ge1\).  Summing over the finitely many
overlap patterns gives
\begin{equation}
        V_2=o(1).
        \label{eq:V2-small}
\end{equation}
Combining \eqref{eq:variance-split}, \eqref{eq:V1-small}, and
\eqref{eq:V2-small}, we obtain
\(\Var(m_k(n))\to0\).
Together with \eqref{eq:expected-moment-limit}, Chebyshev's inequality yields
\begin{equation}
        \frac1n\Tr M_n^k
        \xrightarrow{\Pp}
        \begin{cases}
        0,& k\text{ odd},\\[1mm]
        \Cat_{k/2},& k\text{ even}.
        \end{cases}
        \label{eq:moments-in-prob}
\end{equation}
\subsection{Conclusion by the method of moments}

The limiting moments in \eqref{eq:moments-in-prob} are exactly the moments of
the standard semicircle law:
\[
        \int x^k\,\scm(dx)
        =
        \begin{cases}
        0,& k\text{ odd},\\[1mm]
        \Cat_{k/2},& k\text{ even}.
        \end{cases}
\]
By the standard method-of-moments theorem for random probability measures with a
moment-determinate limiting law \cite{anderson2010introduction}, the empirical
spectral distribution of \(M_n\) converges in probability to \(\scm\).

It remains to transfer this centered convergence to the normalization appearing
in Theorem~\ref{thm:main}.  Put
\(\widetilde M_n=(A-pJ)/\sqrt{npq}\).  Then
\(M_n-\widetilde M_n=pI/\sqrt{npq}\), so \(M_n\) and \(\widetilde M_n\) differ
by a deterministic scalar shift of size
\(\sqrt{p/(nq)}=o(1)\).  Hence \(\mu_{\widetilde M_n}\) also converges in
probability to \(\scm\).  Next,
\(A/\sqrt{npq}-\widetilde M_n=pJ/\sqrt{npq}\), and the right-hand side has rank
one.  By the rank inequality for empirical
spectral distributions \cite{bai2010spectral},
\[
        \sup_x
        \left|
        F_{A/\sqrt{npq}}(x)-F_{\widetilde M_n}(x)
        \right|
        \le \frac1n ,
\]
where \(F_B\) denotes the distribution function of the empirical spectral
distribution of the Hermitian matrix \(B\).  Therefore the empirical spectral
distribution of \(A/\sqrt{npq}\) also converges in probability to \(\scm\).
Finally, \(A/\sqrt{np}=\sqrt q\,A/\sqrt{npq}\), and \(q=1-p\to1\).  Thus the empirical spectral distribution of
\(A/\sqrt{np}\) converges in probability to \(\scm\), proving
Theorem~\ref{thm:main}.

\section{Proof of Theorem~\ref{thm:constantp}}\label{sec:proof_constant}
In the bounded expected degree regime, we do not center the adjacency matrix.
Instead, we expand each edge indicator as \(p+\xi_{ij}\) and compare the resulting
trace moments directly with those of the sparse Erd\H{o}s--R\'enyi graph.
The comparison is especially simple on trees: conditioning from
the leaves shows that their edge indicators have exactly the same joint moment
as independent Bernoulli edges.  Every non-tree support loses a power of \(n\),
and Proposition~\ref{prop:corr} shows that geometric dependence cannot recover
that loss.
Assume
\[
        p=\frac{\alpha}{n},
        \qquad
        \alpha>0\ \text{fixed},
        \qquad
        d=\omega(\log n).
\]
Put \(\eta_{ij}=\1_{\{\ip{X_i}{X_j}\ge\tau\}}\) and
\(\xi_{ij}=\eta_{ij}-p\).  Then \(\eta_{ij}=p+\xi_{ij}\), and
\(L=\log(1/p)=\log n-\log\alpha\).
Since \(d=\omega(L)\), we have \(\eps_n=\sqrt{L/d}=o(1)\), and hence
\begin{equation}
        p^{-D\eps_n}L^D=n^{o(1)}
        \label{eq:p-eps-subpolynomial}
\end{equation}
for every fixed \(D\).
Indeed, \(p^{-D\eps_n}=\exp(D\eps_n L)=\exp(o(L))=n^{o(1)}\), and
\(L^D=n^{o(1)}\).

We first record a consequence of Proposition~\ref{prop:corr} for uncentered
edge indicators.  If \(H\) is a fixed simple graph with \(e\) edges, then
\begin{equation}
        \E\prod_{f\in E(H)}\eta_f
        =
        \sum_{F\subseteq E(H)}
        p^{e-|F|}
        \E\prod_{f\in F}\xi_f .
        \label{eq:uncentered-expansion}
\end{equation}
If \(H\) is a tree, then every nonempty \(F\subseteq E(H)\) is a forest and
therefore has a leaf.  By Proposition~\ref{prop:corr}, all nonempty centered
terms in \eqref{eq:uncentered-expansion} vanish, so
\begin{equation}
        \E\prod_{f\in E(H)}\eta_f=p^e
        \qquad\text{when \(H\) is a tree.}
        \label{eq:tree-exact}
\end{equation}
For a general fixed \(H\), the same expansion and Proposition~\ref{prop:corr}
give constants \(A_H,B_H<\infty\), depending only on \(H\), such that
\begin{equation}
        \left|\E\prod_{f\in E(H)}\eta_f\right|
        \le
        A_H p^{e-B_H\eps_n}L^{B_H}.
        \label{eq:uncentered-bound}
\end{equation}
Indeed, the empty subset contributes \(p^e\); any nonempty subset with a leaf
contributes zero; and each remaining subset is bounded by
\eqref{eq:correlation-bound}.  The finite sum over \(F\subseteq E(H)\) is
absorbed into \(A_H,B_H\).

With this fixed-graph comparison in hand, the trace moments can be compared
term by term with the sparse Erd\H{o}s--R\'enyi moments.

Fix \(k\ge1\).  Expanding the normalized trace gives
\begin{equation}
        \E\frac1n\Tr\left(\frac{A}{\sqrt\alpha}\right)^k
        =
        \frac{1}{n\alpha^{k/2}}
        \sum_w
        \E\prod_{\ell=1}^k\eta_{i_\ell i_{\ell+1}},
        \label{eq:constant-moment-expansion}
\end{equation}
where \(w=(i_1,\ldots,i_k)\) ranges over labeled cyclic words, with
\(i_{k+1}=i_1\) and \(i_\ell\ne i_{\ell+1}\).  Let \(G_w\) be the simple
support graph of \(w\), with \(v(w)\) vertices and \(e(w)\) edges.  Since the
\(\eta_{ij}\)'s are indicators, repeated traversals of the same edge do not
change the product:
\[
        \prod_{\ell=1}^k\eta_{i_\ell i_{\ell+1}}
        =
        \prod_{f\in E(G_w)}\eta_f.
\]
For each fixed unlabeled word pattern, the number of labelings is
\(O_k(n^{v(w)})\).

If \(G_w\) is not a tree, then \(G_w\) is connected and \(e(w)\ge v(w)\).
Using \eqref{eq:uncentered-bound} and \eqref{eq:p-eps-subpolynomial}, the total
contribution of this fixed pattern is bounded by
\begin{equation}
        C_{k,\alpha}\frac{n^{v(w)}p^{e(w)-C_k\eps_n}L^{C_k}}
        {n\alpha^{k/2}}
        \le
        C_{k,\alpha} n^{v(w)-1-e(w)+o(1)}
        =o(1).
        \label{eq:constant-nontree-vanish}
\end{equation}
Here we used \eqref{eq:p-eps-subpolynomial} to absorb both the
correlation loss and the fixed power of \(L\), and then used
\(e(w)\ge v(w)\).
For the Erd\H{o}s--R\'enyi graph, the same non-tree pattern has
expectation \(p^{e(w)}\), so its normalized contribution is
\[
        O_k\left(n^{v(w)-1}p^{e(w)}\right)
        =O_{k,\alpha}\left(n^{v(w)-1-e(w)}\right)
        =O_{k,\alpha}(n^{-1}).
\]
If \(G_w\) is a tree, then \eqref{eq:tree-exact} shows that the contribution is
exactly the same as for the Erd\H{o}s--R\'enyi graph \(\mathcal G(n,\alpha/n)\).
Therefore, for every fixed \(k\),
\begin{equation}
        \E\frac1n\Tr\left(\frac{A}{\sqrt\alpha}\right)^k
        -
        \E\frac1n\Tr\left(\frac{A_{\mathrm{ER}}}{\sqrt\alpha}\right)^k
        \longrightarrow 0.
        \label{eq:constant-expected-comparison}
\end{equation}
For completeness, the convergence of the Erd\H{o}s--R\'enyi trace
moment follows from the same finite-word enumeration.  A non-tree pattern is
\(O_{k,\alpha}(n^{-1})\) by the estimate above.  Let \(\mathcal T_k\) be the
finite set of canonical cyclic word patterns of length \(k\) whose support is a
tree.  If \(P\in\mathcal T_k\) has \(v(P)\) vertices and \(e(P)=v(P)-1\)
edges, then it has \((n)_{v(P)}=(1+o(1))n^{v(P)}\) injective labelings and
contributes
\[
 \frac{(n)_{v(P)}(\alpha/n)^{e(P)}}{n\alpha^{k/2}}
 =\alpha^{e(P)-k/2}+o(1).
\]
Consequently,
\[
 \E\frac1n\Tr\left(\frac{A_{\mathrm{ER}}}{\sqrt\alpha}\right)^k
 \longrightarrow
 m_k(\alpha):=\sum_{P\in\mathcal T_k}\alpha^{e(P)-k/2}.
\]

Jung and Lee~\cite[Lemma~2.2]{jung2018delocalization} prove that these
limiting trace moments are the moments of \(\nu_\alpha\); in particular,
\(m_k(\alpha)=\int x^k\,\nu_\alpha(dx)\).

It remains to show the concentration of these moments.  Let
\(m_k(n)=n^{-1}\Tr(A/\sqrt\alpha)^k\).
For two labeled cyclic words \(w_1,w_2\), let \(V(w_j)\) be the set of labels
visited by \(w_j\).  If \(V(w_1)\cap V(w_2)=\varnothing\), then the two products
of edge indicators depend on disjoint sets of sphere points, so the covariance
is zero.  Hence only pairs sharing at least one vertex need to be considered.
By the triangle inequality,
\begin{equation}
\begin{aligned}
        \Var(m_k(n))
        \le
        \frac{1}{n^2\alpha^k}
        \sum_{\substack{w_1,w_2\\ V(w_1)\cap V(w_2)\ne\varnothing}}
        \biggl(
        \left|\E\prod_{f\in E(G_{w_1}\cup G_{w_2})}\eta_f\right|
        +
        \left|\E\prod_{f\in E(G_{w_1})}\eta_f\right|
        \left|\E\prod_{f\in E(G_{w_2})}\eta_f\right|
        \biggr).
\end{aligned}
        \label{eq:constant-variance-split}
\end{equation}

Fix an overlap pattern for such a pair.  Let \(v_j\) and \(e_j\) be the numbers of
vertices and edges, respectively, in the support of \(w_j\), and let
\(r=|V(w_1)\cap V(w_2)|\ge1\).
The union support has exactly \(v_{12}=v_1+v_2-r\) vertices and
\(e_{12}\) edges.
It is connected, so \(e_{12}\ge v_{12}-1\).  By
\eqref{eq:uncentered-bound}, the contribution of the joint expectation term is
at most
\begin{equation}
\begin{aligned}
        C_{k,\alpha}
        \frac{n^{v_{12}}p^{e_{12}-C_k\eps_n}L^{C_k}}
        {n^2\alpha^k}
        &\le
        C_{k,\alpha} n^{v_{12}-2-e_{12}+o(1)}
        \le
        C_{k,\alpha} n^{-1+o(1)}
        =o(1).
\end{aligned}
        \label{eq:constant-joint-variance}
\end{equation}
For the product of expectations, applying \eqref{eq:uncentered-bound} to the
two supports separately gives
\begin{equation}
\begin{aligned}
        C_{k,\alpha}
        \frac{n^{v_1+v_2-r}
        p^{e_1+e_2-C_k\eps_n}L^{C_k}}
        {n^2\alpha^k}
        &\le
        C_{k,\alpha} n^{v_1+v_2-r-2-e_1-e_2+o(1)}
        \le
        C_{k,\alpha} n^{-r+o(1)}
        =o(1),
\end{aligned}
        \label{eq:constant-product-variance}
\end{equation}
because each support is connected and \(e_j\ge v_j-1\).
Since there are only
finitely many overlap patterns for fixed \(k\), \eqref{eq:constant-joint-variance}
and \eqref{eq:constant-product-variance} imply
\begin{equation}
        \Var(m_k(n))\to0.
        \label{eq:constant-variance}
\end{equation}
Thus the normalized trace moments of \(A/\sqrt{\alpha}\) converge in probability
to the moments of \(\nu_\alpha\).  We verify moment determinacy
directly.  A limiting moment of order \(2r\) is a sum over tree-supported
normalized words of length \(2r\).  Every edge of a tree is a bridge, so a
closed walk traverses it at least twice.  Thus a contributing pattern has at
most \(r\) edges.  There are at most \((2r)^{2r}\) such word patterns, and a
pattern with \(e\le r\) edges has weight
\(\alpha^{e-r}\le C_\alpha^r\), where
\(C_\alpha=\max\{1,\alpha^{-1}\}\).  Hence
\[
        \int |x|^{2r}\,\nu_\alpha(dx)
        \le (2r)^{2r}C_\alpha^r,
\]
and therefore
\[
        \left(\int |x|^{2r}\,\nu_\alpha(dx)\right)^{-1/(2r)}
        \ge \frac{1}{2r\sqrt{C_\alpha}}.
\]
The sum of the right-hand side over \(r\ge1\) diverges, so Carleman's criterion
(see, e.g.,~\cite{bai2010spectral}) shows that \(\nu_\alpha\) is
moment-determinate.
By the method-of-moments theorem, the empirical spectral distribution of
\(A/\sqrt{\alpha}\) converges weakly in probability to \(\nu_\alpha\).  Since
\(np=\alpha\), this is the claimed convergence for \(A/\sqrt{np}\), proving
Theorem~\ref{thm:constantp}.

\paragraph{Acknowledgments}
Y.Z. thanks Manuel Fernandez and Aleksa Milojević for identifying a gap in an earlier version of the manuscript. In developing the proof of Proposition~\ref{prop:corr}, the authors used OpenAI's ChatGPT with the GPT-5.5 Pro model. The authors independently verified the argument and take full responsibility for the final proof. Y.Z. was partially supported by Simons Foundation grant MPS-TSM-00013944. This material is based upon work supported by the Swedish Research Council under grant no. 2021-06594 while Y.Z. was in residence at Institut Mittag-Leffler in Djursholm, Sweden, during the fall of 2024.

\bibliographystyle{plain}
\bibliography{ref}

\appendix
\numberwithin{equation}{section}

\section{One-dimensional spherical estimates}\label{sec:spherical-estimates}
This appendix records the one-dimensional facts about spherical caps used in
Section~\ref{sec:correlation-bound}.

Let \(X,Y\) be independent uniform points on \(S^{d-1}\) and set
\(Z=\sqrt d\,\ip{X}{Y}\) and \(a=\sqrt d\,\tau\).
Then \(\Pp\{Z\ge a\}=p\).  The density of \(Z\) is the standard spherical-coordinate
density \cite{muirhead1982aspects}:
\begin{equation}
        f_d(z)=c_d\left(1-\frac{z^2}{d}\right)^{(d-3)/2}\1_{\{|z|<\sqrt d\}},
        \label{eq:one-density}
\end{equation}
where \(c_d\asymp1\).

Lemma~\ref{lem:tail} collects exactly the one-coordinate estimates used in
Section~\ref{sec:correlation-bound}.  The Gaussian bounds control the
complement of the convergence region, while the relative-density and
shifted-integral estimates handle generators of arbitrary degree.
\begin{appendixlemma}[One-edge density and integral bounds]\label{lem:tail}
Let \(Z=\sqrt d\,\ip{X}{Y}\), \(a=\sqrt d\,\tau\), and
\(L=\log(1/p)\), where \(X,Y\) are independent uniform points on \(S^{d-1}\) and
\(\Pp\{Z\ge a\}=p\).  Assume \(p\to0\) and \(L/d\to0\).  There are absolute
constants \(0<c_*<C_*<\infty\) and \(C>1\) such that, for all sufficiently
large \(n\),
\begin{equation}
        f_d(z)\le C_*e^{-c_*z^2},
        \qquad |z|\le\sqrt d/2,
        \label{eq:basic-gaussian-density}
\end{equation}
and, for every \(t\ge0\),
\begin{equation}
        \Pp\{Z\ge t\}\le C_*e^{-c_*t^2}.
        \label{eq:basic-gaussian-tail}
\end{equation}
The cap threshold satisfies
\begin{equation}
        0\le a\le C_*\sqrt L,
        \label{eq:a-bound}
\end{equation}
and its density obeys the relative estimate
\begin{equation}
        f_d(z)\le C_*(1+a)p\,e^{-c_*(z^2-a^2)},
        \qquad a\le z<\sqrt d.
        \label{eq:density-decay-from-a}
\end{equation}
For every integer \(r\ge0\),
\begin{equation}
        (1+a)\int_a^\infty z^r e^{-c_*(z^2-a^2)}\,dz
        \le C^{r+1}(a+\sqrt r)^r.
        \label{eq:shifted-gaussian-integral}
\end{equation}
\end{appendixlemma}

\begin{proof}
By \eqref{eq:one-density},
\(c_d=\Gamma(d/2)/(\sqrt{\pi d}\,\Gamma((d-1)/2))\).
By Stirling's formula, there are absolute constants
\(0<c_*<C_*<\infty\) such that \(c_*\le c_d\le C_*\) for all \(d\ge2\).

We first derive the Gaussian bounds.  If \(|z|\le\sqrt d/2\), then
\(\log(1-z^2/d)\le -z^2/d\), and therefore
\[
        f_d(z)\le C_* e^{-c_*z^2},
        \qquad |z|\le\sqrt d/2,
\]
which proves \eqref{eq:basic-gaussian-density}.
If \(z\ge\sqrt d/2\), then
\((1-z^2/d)^{(d-3)/2}\le C_*e^{-c_*d}\).  Integrating the two
density bounds gives the tail estimate explicitly as follows.
If \(0\le t\le\sqrt d/2\), then, after decreasing the positive
constant in the exponent,
\[
\begin{aligned}
 \Pp\{Z\ge t\}
 &\le C\int_t^\infty e^{-cz^2}\,dz
      +C\sqrt d\,e^{-cd} \\
 &\le C_*e^{-c_*t^2}.
\end{aligned}
\]
If \(\sqrt d/2<t<\sqrt d\), then
\(\Pp\{Z\ge t\}\le C\sqrt d\,e^{-cd}\le C_*e^{-c_*t^2}\); if
\(t\ge\sqrt d\), the probability is zero.
Thus, for every \(t\ge0\),
\[
        \Pp\{Z\ge t\}\le C_* e^{-c_*t^2}.
\]
This proves \eqref{eq:basic-gaussian-tail}.
Since \(p\to0\), we have \(p<1/2\) for all sufficiently large \(n\).  The density is symmetric, so the upper
\(p\)-quantile \(a\) is nonnegative.  Using
\(p=e^{-L}=\Pp\{Z\ge a\}\) in \eqref{eq:basic-gaussian-tail} gives
\(a^2\le C_*(L+1)\).  Since \(L\to\infty\), enlarging the absolute constant
gives \(a^2\le C_*L\), proving
\eqref{eq:a-bound}.

We next prove the relative density estimate
\eqref{eq:density-decay-from-a}.  Since \(a\le C_*\sqrt L\) and \(L/d\to0\),
we have \(a\le \sqrt d/4\) for all sufficiently large \(n\).  The logarithmic derivative of \(f_d\) is
\begin{equation}
        \frac{d}{dz}\log f_d(z)
        =-\frac{(d-3)z}{d-z^2}.
        \label{eq:log-derivative-density}
\end{equation}
Set \(I_a=[a,a+(1+a)^{-1}]\).  If \(z\in I_a\), then
\(z\le a+1\).  Since \(a\le\sqrt d/4\) and \(d\to\infty\), we have
\(z\le\sqrt d/2\) for all sufficiently large \(n\), and hence
\(d-z^2\ge3d/4\).  It follows from
\eqref{eq:log-derivative-density} that
\[
\left|\frac{d}{dz}\log f_d(z)\right|
=\frac{(d-3)z}{d-z^2}
\le \frac{dz}{3d/4}
\le \frac43(1+a).
\]
Consequently, for every \(z\in I_a\),
\[
\log f_d(z)-\log f_d(a)
\ge -\frac43(1+a)(z-a)
\ge -\frac43.
\]
Thus \(f_d(z)\ge e^{-4/3}f_d(a)\) on \(I_a\).  After decreasing the absolute
constant \(c_*\), if necessary, we obtain
\[
        p=\int_a^{\sqrt d}f_d(z)\,dz
        \ge \int_{I_a}f_d(z)\,dz
        \ge \frac{c_* f_d(a)}{1+a},
\]
and so
\[
        f_d(a)\le C_*(1+a)p.
\]
For \(d\ge6\) and \(a\le z<\sqrt d\),
\[
        \frac{d-3}{d-z^2}\ge\frac{d-3}{d}\ge\frac12.
\]
Therefore
\[
        \frac{d}{dz}\log f_d(z)
        =-\frac{(d-3)z}{d-z^2}\le-\frac z2.
\]
Integrating from \(a\) to \(z\) gives
\[
        f_d(z)
        \le f_d(a)e^{-c_*(z^2-a^2)}
        \le C_*(1+a)p\,e^{-c_*(z^2-a^2)}.
\]
This is \eqref{eq:density-decay-from-a}.

It remains to prove the shifted Gaussian integral estimate.  For \(r=0\),
set \(u=z-a\).  If \(0\le a\le1\), then
\[
 \int_0^\infty e^{-c_*(u^2+2au)}\,du
 \le\int_0^\infty e^{-c_*u^2}\,du
 \le\frac{C}{1+a}.
\]
If \(a>1\), then
\[
 \int_0^\infty e^{-c_*(u^2+2au)}\,du
 \le\int_0^\infty e^{-2c_*au}\,du
 =\frac{1}{2c_*a}
 \le\frac{C}{1+a}.
\]
Thus
\[
        \int_0^\infty e^{-c_*(u^2+2au)}\,du
        \le \frac{C}{1+a}.
\]
For \(r\ge1\), put
\[
        J_r(a)=\int_0^\infty u^r e^{-c_*(u^2+2au)}\,du.
\]
Using \((a+u)^r\le2^{r-1}(a^r+u^r)\), we obtain
\[
\begin{aligned}
        (1+a)\int_a^\infty z^r e^{-c_*(z^2-a^2)}\,dz
        &\le
        2^{r-1}(1+a)
        \left(
        a^r\frac{C}{1+a}+J_r(a)
        \right).
\end{aligned}
\]
If \(a\le\sqrt r\), then
\[
        J_r(a)
        \le \int_0^\infty u^r e^{-c_*u^2}\,du
        =\frac12c_*^{-(r+1)/2}\Gamma\left(\frac{r+1}{2}\right)
        \le C^{r+1}r^{r/2}.
\]
The first term in the right-hand side of the preceding bound is at
most \(C^ra^r\).  For the second term, \(a\le\sqrt r\) and
\(1+\sqrt r\le C^r\) for \(r\ge1\), so
\[
        2^{r-1}(1+a)J_r(a)
        \le C^{r+1}r^{r/2}
        \le C^{r+1}(a+\sqrt r)^r.
\]
This proves \eqref{eq:shifted-gaussian-integral} when \(a\le\sqrt r\).
If \(a>\sqrt r\), then
\[
        J_r(a)
        \le \int_0^\infty u^r e^{-2c_*au}\,du
        =\frac{r!}{(2c_*a)^{r+1}}.
\]
Since \(a>\sqrt r\ge1\), we have \((1+a)/a\le2\) and
\(r!\le r^r\).  Hence
\[
        (1+a)J_r(a)
        \le C^{r+1}\left(\frac ra\right)^r
        \le C^{r+1}a^r,
\]
where the last inequality uses \(r<a^2\).  Together with the \(a^r\) term in
the preceding decomposition, this proves
\eqref{eq:shifted-gaussian-integral} in the second case.
This completes the proof.
\end{proof}

\end{document}